\documentclass[3p,numbers,sort&compress, nopreprintline]{elsarticle}
\nonstopmode
\usepackage{lineno,hyperref,amsmath,amssymb, physics, filecontents, xcolor, relsize}

\usepackage[nodocdata = true, noproducerdata=true, noeditdata = true, noptexdata=true, nopdftrailerid=true]{pdfprivacy}

\journal{}

\usepackage{mathrsfs}
\usepackage{setspace}
\usepackage[]{algorithm}
\usepackage[]{algpseudocode}

\renewcommand{\vec}[1]{\ensuremath{\boldsymbol{#1}}}

\newcommand{\vecstate}[1]{\underline{\vec{#1}}}
\newcommand{\scastate}[1]{\underline{#1}}
\newcommand{\angleb}[1]{\langle #1 \rangle}
\newcommand{\squareb}[1]{\left[#1\right]}

\usepackage{numcompress}

\date{}

\begin{document}

\begin{frontmatter}

\title{A stabilized computational nonlocal poromechanics model for dynamic analysis of saturated porous media}

\author{Shashank Menon}
\author{Xiaoyu Song \corref{mycorrespondingauthor}}
\address{Engineering School of Sustainable Infrastructure and Environment\\
 University of Florida, Gainesville, FL 32611 USA}

\cortext[mycorrespondingauthor]{Corresponding author}
\ead{xysong@ufl.edu}

\begin{abstract}

In this article we formulate a stable computational nonlocal poromechanics model for dynamic analysis of saturated porous media. As a novelty, the stabilization formulation eliminates zero-energy modes associated with the original multiphase correspondence constitutive models in the coupled nonlocal poromechanics model. The two-phase stabilization scheme is formulated based on an energy method that incorporates inhomogeneous solid deformation and fluid flow. In this method, the nonlocal formulations of skeleton strain energy and fluid flow dissipation energy equate to their local formulations. The stable coupled nonlocal poromechanics model is solved for dynamic analysis by an implicit time integration scheme. As a new contribution, we validate the coupled stabilization formulation by comparing numerical results with analytical and finite element solutions for one-dimensional and two-dimensional dynamic problems in saturated porous media. Numerical examples of dynamic strain localization in saturated porous media are presented to demonstrate the efficacy of the stable coupled poromechanics framework for localized failure under dynamic loads.

\end{abstract}

\begin{keyword}
stabilization \sep nonlocal \sep  coupled \sep poromechanics \sep dynamics \sep saturated porous media 
\end{keyword}
\end{frontmatter}
%

\section{Introduction}

Dynamic behavior of deforming porous media is significant problem in engineering and science \cite{zienkiewicz1999computational, lewis1998finite}. A fully coupled dynamic analysis of porous media is essential in earthquake engineering (soil liquefaction), geo-hazard engineering (dynamics landslide triggering), and offshore wind industry (wind turbine foundations) \cite{Zienkiewicz1990Static, alonso2021triggering}. The integrity of civil infrastructure (e.g., dams and levee systems) can be seriously compromised by the poor performance of soils under dynamic loading conditions. Thus, numerous researchers have studied the dynamics and wave propagation characteristics of saturated porous media, e.g., \cite{Prevost1985Wave, vardoulakis1986dynamic, diebels1996dynamic, alonso2003influence, popescu2006dynamics} among others. One typical failure of porous media under dynamic loads is dynamic strain localization (e.g., \cite{Loret1991Dynamic, lewis1998finite}). A topic has been insufficiently studied compared to the extensive research on strain localization of porous media under static or quasi-static loading conditions (see \cite{song2017strain} and \cite{wang2020strain} for a recent review). Over the past decades, viscoplasticity has been adopted to study dynamic strain localization of solid or porous media (see \cite{needleman1988material, Loret1991Dynamic, schrefler1996multiphase, heider2014dynamic, shahbodagh2014dynamic, oka2019computational}, among others), as well as a means of regularizing the rate-independent problem so that the governing equation of dynamic problems remains hyperbolic. Heider et al \cite{heider2014dynamic} developed a a coupled dynamic elasto-viscoplastic model for sands and found that inertial loads had a significant impact on dynamic strain localization. Shahbodagh et al \cite{shahbodagh2014dynamic} proposed an elasto-viscoplastic model for dynamic analysis of strain localization in fully saturated clay. The above methods are based on the classic local poromechanics. It has been generally recognized that nonlocal models are robust to investigate mulitphysics failure mechanisms including the dynamic strain localization phenomenon in porous media (e.g., \cite{cosserat1909theorie, eringen1964nonlinear, kroner1967elasticity, DeBorst1991Simulation} and many others). In this article, we propose a stabilized coupled nonlocal poromechanics model to study dynamic strain localization in saturated porous media.

In \cite{menon2020computational}, a computational periporomechanmics model was formulated to model localized failure in unsaturated porous media under static condition. The coupled nonlocal model in \cite{menon2020computational} was formulated based on the peridynamic state concept \cite{Silling2007Peridynamic}, the effective force state concept and multiphase correspondence principle \cite{song2020peridynamic}. We refer to the literature for other nonlocal models for porous media (see \cite{turner2013non, madenci2014peridynamic, jabakhanji2015peridynamic, ouchi2015fully, Oterkus2017Fully, menon2019coupled, zhang2019coupling, ni2020hybrid}, among others), which are formulated using peridynamics theory (i.e., the bond-based or ordinary state-based peridynamics) \cite{Silling2000Reformulation, Silling2007Peridynamic, madenci2014peridynamic} and poroelasticity (see \cite{cheng2016poroelasticity}). It was demonstrated that the coupled nonlocal model in \cite{song2020peridynamic} is robust for modeling localized failure in unsaturated soils under static/quasi-static loading conditions. However, it is recognized that the peridynamics model formulated based on the correspondence principle for the single-phase solid has stability issues under extreme large deformation and dynamic loading conditions (e.g., \cite{littlewood2011nonlocal, breitenfeld2014non, extendedTupek2014, silling2017stability, bobaru2016handbook, Li2018stabilized, Gu2018Revisit, Wu2015stabilized, Chen2018Bond,RoyChowdhury2019modified, hashim2020implicit}, among others). Silling \cite{silling2017stability} showed that the numerical oscillation in peridynamics material models formulated via the single-phase correspondence principle \cite{Silling2007Peridynamic} is associated with material instability instead of merely an artifact of the meshless discretization \cite{silling2005meshfree}. This study suggests that the primary causes of instability and numerical oscillations are (1) the weak dependence of the force density in a bond on its own deformation and (2) the loss of the non-uniform part of deformation due to the integration over the horizon of a material point. Note that the above studies and remedies for the instability were focused on modeling the large deformation or extreme damage of single-phase solid materials through the original correspondence principle \cite{Silling2007Peridynamic}. Here for the first time we prove that the recently formulated multiphase correspondence for modeling unsaturated porous media inherit the aforementioned instability under large deformation and dynamic loading conditions. To circumvent the issue, we have formulated the coupled stabilization terms for solid deformation and fluid flow respectively in Section 2.

In this article, we formulate and implement a stabilized coupled nonlocal poromechanics model for dynamic strain localization in saturated porous media. As a novelty, the stabilization terms eliminates zero-energy modes associated with the multiphase correspondence constitutive models in the coupled nonlocal framework for modeling saturated porous media under dynamic loads. The stabilization terms are formulated based on the energy method in which the nonlocal formulations of skeleton strain energy and fluid flow dissipation energy equate to their local formulations in line with the classical poromechanics for saturated porous media. Specifically, the method incorporates non-homogeneous solid deformation and fluid flow around a material point. The stabilized coupled nonlocal poromechanics model is numerically solved using a hybrid Lagrangian-Eulerian meshless method with an implicit time integration scheme. Parallel computing is also adopted for computational efficiency. The coupled stabilized formulation is validated by comparing numerical results with analytical and finite element solutions for one-dimensional and two-dimensional dynamic problems in saturated porous media. Numerical examples of dynamic strain localization in saturated porous media are presented to demonstrate the efficacy of the stable coupled poromechanics framework for localized failure under dynamic loads. We note that the stabilized coupled nonlocal model for saturated porous media can be readily extended to model dynamic problems in unsaturated porous media. 

The contribution of this article includes (1) a proof of the zero-energy modes associated with the original multiphase correspondence constitutive principle for modeling porous media through multiphase peridynamic states, (2) a remedy based on the energy method to remove the multiphase zero-energy modes for dynamic analysis, and (3) an implicit numerical implementation of the proposed stabilized dynamic non-local poromechanics model and its validation. For sign convention, the assumption in continuum mechanics is followed, i.e., for solid skeleton, tensile force/stress is positive and compression is negative,  and for fluid compression is positive and tension is negative.

\section{Stabilized nonlocal formulation for  coupled dynamics problems}

For conciseness of notations, it is assumed that the peridynamic state variable without a prime means that the variable is evaluated at $\vec{x}$ on the associated bond $\vec{x}' - \vec{x}$ and the peridynamic state variable with a prime means that the variable is evaluated at $\vec{x}'$ on the associated bond $\vec{x} - \vec{x}'$, e.g., $\vecstate{T} = \vecstate{T}[\vec{x}]\langle\vec{x}'-\vec{x}\rangle$ and $\vecstate{T}' = \vecstate{T}[\vec{x}']\langle\vec{x}-\vec{x}'\rangle$.

\subsection{Dynamic saturated periporomechanics model}
\begin{figure}[h]
\centering
  \includegraphics[width=0.5\textwidth]{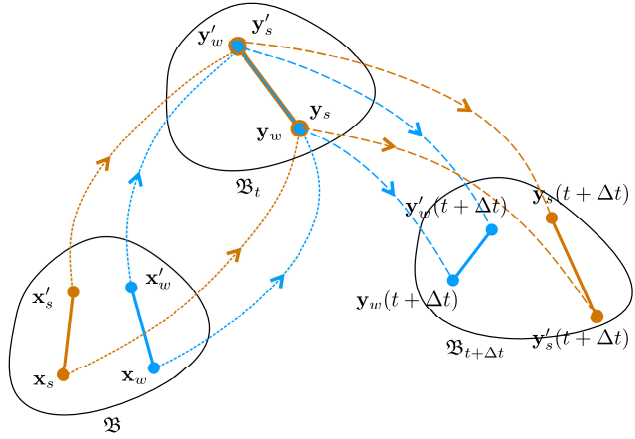}
\caption{Schematics of kinematics of two material points $\vec{x}$ and $\vec{x}'$.}
\label{kinematics_peri-poromechanics}   
\end{figure}

Periporomechanics is a fully coupled, nonlocal theory of porous media. It is a reformulation of classical poromechanics through peridynamics for modeling continuous or discontinous deformation and physical processes in porous media \cite{song2020peridynamic, menon2020computational}. In periporomechanics, it is assumed that a porous media body is composed of material points which have two kinds of degree of freedom, i.e., displacement and fluid pressure. A material point $\vec{x}$ has poromechanical and physical interactions with any material point $\vec{x}'$ within its neighborhood, $\mathcal{H}$. Here $\mathcal{H}$ is a spherical domain around $\vec{x}$ with radius $\delta$, which is called the horizon for the porous medium. Let $\rho$ be the density of the two phase mixture that is determined by
\begin{equation}
\rho = \rho_s (1-\phi) + \rho_w \phi,
\label{mixture_density}
\end{equation}
where $\rho_s$ and $\rho_w$ are the intrinsic density of the solid and the fluid, respectively, and $\phi$ is the porosity (i.e., the volume of pore space divided by the total volume). Figure \ref{kinematics_peri-poromechanics} shows the schematics of the kinematics of two material points. 

The relative position of material points $\vec{x}$ and $\vec{x}'$ in the reference configuration is denoted by
\begin{equation}
\vec{\xi} = \vec{x}' - \vec{x}.
\end{equation}
Let the displacements of $\vec{x}$ and $\vec{x}'$ be $\vec{u}$ and $\vec{u}'$ respectively. The relative displacement between $\vec{x}$ and $\vec{x}'$ in the deformed configuration is
\begin{equation}
\vec{\eta} = \vec{u}' - \vec{u}.
\end{equation}
The position vectors $\vec{y}'$ and $\vec{y}$ of $\vec{x}'$ and $\vec{x}$ in the deformed configuration are
\begin{align}
\vec{y} &= \vec{x}  + \vec{u},\quad \text{and}\quad \vec{y}' = \vec{x}'  + \vec{u}'.
\end{align}
The deformation state of $\vec{x}$ on $\vec{\xi}$ is defined as
\begin{equation}
\vecstate{Y} = \vec{y}' - \vec{y} = \vec{\xi} + \vec{\eta}.
\end{equation}
The fluid potential state at $\vec{x}$ imposed on $\vec{\xi}$ is defined as
\begin{equation}
\scastate{\Phi}\angleb{\vec{\xi}} = p'_w - p_w,
\end{equation}
where $p'_w$ and $p_w$ are pore fluid pressures at $\vec{x}'$ and $\vec{x}$, respectively.

\begin{equation}
\int_{\mathcal{H}}\left\{ \left( \overline{\vecstate{T}}\squareb{\vec{x}}\angleb{\vec{\xi}}- \vecstate{\vec{T}}_w\squareb{\vec{x}}\angleb{\vec{\xi}} \right) - \left( \overline{\vecstate{T}}\squareb{\vec{x}'}\angleb{-\vec{\xi}}- \vecstate{\vec{T}}_w\squareb{\vec{x}'}\angleb{-\vec{\xi}} \right) \right\} dV' + \rho \vec{g} 
= \rho \ddot{\vec{u}},
\label{momentum_balance_total_force_state}
\end{equation} 
where $\overline{\vecstate{T}}\squareb{\vec{x}}\angleb{\vec{\xi}}$ and $ \vecstate{\vec{T}}_w\squareb{\vec{x}}\angleb{\vec{\xi}}$ the effective state and the fluid phase force state at material point $\vec{x}$, $\vec{g}$ is the gravity acceleration, and $\ddot{\vec{u}}$ is the acceleration vector. Assuming an incompressible solid phase, the fluid mass balance equation incorporating the volume change of the solid skeleton reads
\begin{equation}
\frac{\phi}{K_w}\dot{p}_w + \int_{\mathcal{H}} \left(\dot{\scastate{\mathcal{V}}}\squareb{\vec{x}}\angleb{\vec{\xi}} - \dot{\scastate{\mathcal{V}}}'\squareb{\vec{x}'}\angleb{\vec{-\xi}}\right)dV' + \frac{1}{\rho_w}\int_{\mathcal{H}} \left(\scastate{{\mathcal{Q}}}\squareb{\vec{x}}\angleb{\vec{\xi}}- \scastate{{\mathcal{Q}}}'\squareb{\vec{x}'}\angleb{-\vec{\xi}}\right)dV' = 0
\label{mass_balance_general}
\end{equation}
where $p_w$ is fluid pressure, $K_w$ is the bulk modulus of fluid, $\dot{\scastate{\mathcal{V}}}_s$ and $\dot{\scastate{\mathcal{V}}}'_s$ are the rate of volume change scalar states of the solid at $\vec{x}$ and $\vec{x}'$, respectively, $\scastate{{\mathcal{Q}}}$ and $\scastate{{\mathcal{Q}}}'$ are fluid mass flow scalar states relative to the solid skeleton at $\vec{x}$ and $\vec{x}'$, respectively. 

Through the multiphase correspondence principle \cite{song2020peridynamic}, the effective force state and the fluid mass flow state can be determined by the peridynamic constitutive models via the deformation state of the solid skeleton and the fluid potential state, respectively.
The momentum balance equation for saturated porous media in periporomechanics under dynamic loading reads
\begin{align}
\int_\mathcal{H}\left[\left(\scastate{\omega}(\overline{\vec{\sigma}} - p_w{\vec{1}})J\widetilde{\vec{F}}^{-T}\vec{K}^{-1}\vec{\xi}\right) -\left(\scastate{\omega}(\overline{\vec{\sigma}}' - p'_w{\vec{1}})J\widetilde{\vec{F}}^{-T}\vec{K}^{-1}\vec{\xi}'\right)\right] \text{d}V' + \rho\vec{g}= \rho\ddot{\vec{u}},
\label{momentum_balance_correspondence}
\end{align}
where $\scastate{\omega}$ is a scalar influence function, $\widetilde{\vec{F}}$ is the nonlocal deformation gradient, $J$ is the determinant of $\widetilde{\vec{F}}$, $\overline{\vec{\sigma}}$ is the effective stress tensor that can be determined by classical constitutive models given $\widetilde{\vec{F}}$, $\vec{1}$ is the second-order identity tensor, and $\vec{K}$ is the shape tensor. The nonlocal deformation gradient $\widetilde{\vec{F}}$ is defined as
\begin{equation}
\widetilde{\vec{F}} = \left( \int_{\mathcal{H}}  \scastate{\omega}\vecstate{Y}\otimes \vec{\xi} \text{d}V' \right) \vec{K} ^{-1},\;\; \text{with}\; \vec{K} = \int_{\mathcal{H}} \scastate{\omega}\vec{\xi}\otimes\vec{\xi}\text{d}V '.
\label{nonlocal_deformation_gradient}
\end{equation}

It is assumed that the accelerations of the fluid and the solid are identical \cite{zienkiewicz1999computational}. Through the multiphase correspondence principle \cite{song2020peridynamic}, the fluid mass balance equation accounting for the volume change of the solid skeleton reads
\begin{align}
\frac{\phi}{K_w}\dot{p}_w + \int_{\mathcal H} \left(\scastate{\omega}\dot{\vec{u}} \vec{K}^{-1}\vec{\xi} - \scastate{\omega}\dot{\vec{u}}' \vec{K}^{-1}\vec{\xi}'\right)dV' 
+ \int_{\mathcal{H}} \left[\scastate{\omega} (\vec{q}_w - k_w\ddot{\vec{u}})\vec{K}^{-1}\vec{\xi}   -  \scastate{\omega} (\vec{q}'_w - k_w\ddot{\vec{u}}')\vec{K}^{-1}\vec{\xi} '\right]\text{d}V'= 0,
\label{mass_balance_correspondence}
\end{align}
where $K_w$ is the bulk modulus of fluid, $\vec{q}_w$ and $\vec{q}'_w$ are the fluid flux vectors at $\vec{x}$ and $\vec{x}'$ respectively and $k_w$ is the hydraulic conductivity.
The fluid flux can be determined by Darcy's law as
\begin{equation}
\vec{q}_w = - k_w \widetilde{\vec{\nabla}\Phi},
\label{nonlocal_fluid_flux_uniform}
\end{equation}
where $k_w$ is hydraulic conductivity, and $\widetilde{\vec{\nabla}\Phi}$ is the approximate nonlocal fluid pressure gradient,
\begin{equation}
\widetilde{\vec{\nabla}\Phi} = \left(\int_{\mathcal{H}} \scastate{\omega}  \scastate{\Phi}\vec{\xi}\text{d} V'\right)\vec{K}^{-1}.
\label{nonlocal_pressure_gradient}
\end{equation}

\subsection{Two-phase stabilization formulation}

As stated in the multiphase correspondence principle \cite{song2020peridynamic}, the nonlocal deformation gradient and fluid pressure gradient are approximate for the non-uniform deformation and fluid flow while they are exact for the uniform deformation and fluid flow respectively. There can be increments $\text{d}\vecstate{Y}$ that may have no effect on the approximate deformation gradient, leading to zero-energy modes of deformation \cite{silling2017stability}. Similarly, for fluid flow there can be increments $\text{d}\scastate{\Phi}_w$ that may have no effect on the approximate fluid pressure gradient, leading to zero-energy modes of fluid flow. It can be demonstrated as follows. The nonuniform part of the solid deformation state  and fluid potential state can be defined as 
\begin{align}
\vecstate{\mathcal{R}}^s\angleb{\vec{\xi}} &= \vecstate{Y}\angleb{\vec{\xi}} - \widetilde{\vec{F}}\vec{\xi}
\label{nonuniform_deformation_state},\\
\scastate{\mathcal{R}}^w\angleb{\vec{\xi}} & = \scastate{\Phi}_w \angleb{\vec{\xi}} - \widetilde{\vec{\nabla}\Phi}_w\vec{\xi}
\label{nonunifrom_fluid_potential_state}
\end{align}
Substituting equations \eqref{nonuniform_deformation_state} and \eqref{nonunifrom_fluid_potential_state} into \eqref{nonlocal_deformation_gradient} and \eqref{nonlocal_pressure_gradient} generates

\begin{align}
\left(\int_{\mathcal{H}} \scastate{\omega}\vecstate{\mathcal{R}}^s\otimes\vec{\xi} \text{d}V'\right)\vec{K}^{-1}
&= \left(\int_\mathcal{H}\scastate{\omega}\left(\vecstate{Y} - \widetilde{\vec{F}}\vec{\xi}\right)\otimes\vec{\xi}\text{d}V'\right) \vec{K}^{-1}\nonumber\\
& = \left(\int_\mathcal{H}\scastate{\omega}\vecstate{Y} \otimes\vec{\xi}\text{d}V' \right)\vec{K}^{-1} - \widetilde{\vec{F}} \left(\int_\mathcal{H}\scastate{\omega}\vec{\xi}\otimes\vec{\xi}\text{d}V' \right)\vec{K}^{-1}\nonumber\\
&= \widetilde{\vec{F}} - \widetilde{\vec{F}}\vec{K}\vec{K}^{-1} \nonumber\\
& = \vec{0},
\label{non_uniform_deformation_gradient}
\end{align}

\begin{align}
\left(\int_{\mathcal{H}} \scastate{\omega}\scastate{\mathcal{R}}^w \vec{\xi} \text{d}V'\right)\vec{K}^{-1}
&= \left(\int_\mathcal{H}\scastate{\omega}\left(\scastate{\Phi}_w - \widetilde{\vec{\nabla}\Phi}_w\vec{\xi}\right)\vec{\xi}\text{d}V'\right) \vec{K}^{-1}\nonumber\\
& = \left(\int_\mathcal{H}\scastate{\omega}\scastate{\Phi}\vec{\xi}\text{d}V' \right)\vec{K}^{-1} - \widetilde{\vec{\nabla}\Phi}_w \left(\int_\mathcal{H}\scastate{\omega}\vec{\xi}\otimes\vec{\xi}\text{d}V' \right)\vec{K}^{-1}\nonumber\\
&= \widetilde{\vec{\nabla}\Phi}_w - \widetilde{\vec{\nabla}\Phi}_w\vec{K}\vec{K}^{-1} \nonumber\\
& = 0.
\label{non_uniform_fluid_potential_gradient}
\end{align}

For the solid skeleton, the strain energy density is assumed as
\begin{equation}
\mathscr{W} = \mathscr{W}_c + \mathscr{W}_s,
\label{solid_stored_energy}
\end{equation}
where $\mathscr{W}_c$ is the stored energy from the standard correspondence material model, and  
\begin{equation}
\mathscr{W}_s(\vecstate{Y}) = \dfrac{1}{2}(\scastate{\beta}\vecstate{\mathcal{R}}^s)\bullet\vecstate{\mathcal{R}}^s,
\label{solid_stored_energy_non_uniform_deformation}
\end{equation}
where $\scastate{\beta}$ is a prescribed positive-valued scalar state. Given a small increment $\text{d}\vecstate{Y}$, it follows from \eqref{nonuniform_deformation_state} and \eqref{solid_stored_energy_non_uniform_deformation} that 
\begin{align}
\text{d}\mathscr{W}_s  &= \int_{\mathcal{H}}\scastate{\beta}\scastate{\mathcal{R}}^s_i\left(\text{d}\scastate{Y}_i - \text{d}\widetilde{F}_{ij}\xi_j\right)\text{d}V' \nonumber \\
 &= \int_{\mathcal{H}}\scastate{\beta}\scastate{\mathcal{R}}^s_i\text{d}\scastate{Y}_i\text{d}V' -\int_{\mathcal{H}}\scastate{\beta}\scastate{\mathcal{R}}^s_i \text{d}\widetilde{F}_{ij}\xi_j\text{d}V'\nonumber\\
 & = \int_{\mathcal{H}}\scastate{\beta}\scastate{\mathcal{R}}^s_i\text{d}\scastate{Y}_i\text{d}V' -\left(\int_{\mathcal{H}}\scastate{\beta}\scastate{\mathcal{R}}^s_i \xi_j\text{d}V'\right)\left(\int_{\mathcal{H}}\scastate{\omega}\text{d}Y_i\xi_l\text{d}V'\right)K^{-1}_{lj}\nonumber\\
  & = \int_{\mathcal{H}}\left[\scastate{\beta}\scastate{\mathcal{R}}^s_i-\left(\int_{\mathcal{H}}\scastate{\beta}\scastate{\mathcal{R}}^s_i \xi_j\text{d}V'\right)\scastate{\omega}\xi_l K^{-1}_{lj}\right]\text{d}\scastate{Y}_i\text{d}V',
  \label{energy_solid_non_uniform}
\end{align}
where $i, j, l =1,2,3$. It follows from \eqref{energy_solid_non_uniform} that the stablized term for the effective force state accounting for the non-uniformed deformation state can be written as
\begin{equation}
\overline{\scastate{T}}^s_i = \scastate{\beta}\scastate{\mathcal{R}}^s_i-\left(\int_{\mathcal{H}}\scastate{\beta}\scastate{\mathcal{R}}^s_i \xi_j\text{d}V'\right)\scastate{\omega}\xi_l K^{-1}_{lj}.
\label{effective_force_stabilized_term}
\end{equation}
Given $\text{d}\vecstate{Y}$ it can be proved that $\text{d}\mathscr{W}_s$ is always positive for a positive value of $\scastate{\beta}$. It follows from \eqref{effective_force_stabilized_term} along with the original multiphase correspondence principle \cite{song2020peridynamic} that the stabilized effective stress state can be written as
\begin{equation}
\overline{\vecstate{T}} = \scastate{\omega} J\overline{\vec{\sigma}}\widetilde{\vec{F}}^{-T} \vec{K}^{-1} \vec{\xi} +  \scastate{\beta}\vecstate{\mathcal{R}}^s - \left(\int_{\mathcal{H}}\scastate{\beta}\vecstate{\mathcal{R}}^s\otimes\vec{\xi}\text{d}V'\right)\scastate{\omega}\vec{K}^{-1}\vec{\xi}.
\label{effective_force_state_stabilization_full_form}
\end{equation}
To simplify \eqref{effective_force_state_stabilization_full_form}, it is assumed that
\begin{equation}
\scastate{\beta} = \dfrac{GC}{\omega_0}\scastate{\omega},\quad\text{with}\; \omega_0 = \int_{\mathcal{H}}\scastate{\omega}\text{d}V',
\label{ad_hoc_beta}
\end{equation}
where $G$ is a positive constant on the order of 1 and $C$ is a micromodulus which will be provided in the next section. With \eqref{ad_hoc_beta} and \eqref{non_uniform_deformation_gradient}, it can be proved that the third term of \eqref{effective_force_state_stabilization_full_form} becomes zero. Thus, equation \eqref{effective_force_state_stabilization_full_form} can be expressed as
\begin{equation}
\overline{\vecstate{T}} = \scastate{\omega} \left(J\overline{\vec{\sigma}}\widetilde{\vec{F}}^{-T} \vec{K}^{-1} \vec{\xi} +  \dfrac{GC}{\omega_0}\vecstate{\mathcal{R}}^s\right).
\label{effective_force_state_stabilization_simple_form}
\end{equation}

Similarly, for the fluid phase the energy dissipation accounting for the non-uniform fluid potential state can be written as
\begin{equation}
\mathscr{V}= \mathscr{V}_c + \mathscr{V}_s,
\label{fluid_dissipation_energy_total}
\end{equation}
where $\mathscr{V}_c$ is the energy dissipation related to the uniform fluid potential through the corresponding Darcy's law (i.e., \eqref{nonlocal_fluid_flux_uniform}) and $\mathscr{V}_s$ is the dissipation energy through the non-uniform fluid potential. It is assumed that $\mathscr{V}_s$ takes the general form
\begin{equation}
\mathscr{V}_s = \dfrac{1}{2}\left(\scastate{\lambda}\scastate{\mathcal{R}}^w\right)\cdot\scastate{\mathcal{R}}^w,
\label{fluid_dissipation_energy_non_uniform}
\end{equation}
where $\scastate{\lambda}$ is a positive-valued scalar state. Given a small increment $\text{d}\scastate{\Phi}$, it follows from \eqref{non_uniform_fluid_potential_gradient} and \eqref{fluid_dissipation_energy_non_uniform} that 
\begin{align}
\text{d}\mathscr{V}_s  &= \int_{\mathcal{H}}\scastate{\lambda}\scastate{\mathcal{R}}^w\left(\text{d}\scastate{\Phi} - \text{d}(\widetilde{\vec{\nabla}\Phi})_i\xi_i\right)\text{d}V' \nonumber \\
 &= \int_{\mathcal{H}}\scastate{\lambda}\scastate{\mathcal{R}}^w \text{d}\scastate{\Phi}\text{d}V' -\int_{\mathcal{H}}\scastate{\lambda}\scastate{\mathcal{R}}^w \text{d}(\widetilde{\vec{\nabla}\Phi})_i\xi_i\text{d}V'\nonumber\\
 & = \int_{\mathcal{H}}\scastate{\lambda}\scastate{\mathcal{R}}^w \text{d}\scastate{\Phi}\text{d}V'-\left(\int_{\mathcal{H}}\scastate{\lambda}\scastate{\mathcal{R}} \xi_i\text{d}V'\right)\left(\int_{\mathcal{H}}\scastate{\omega}\text{d}\scastate{\Phi}\xi_j\text{d}V'\right)K^{-1}_{ji}\nonumber\\
  & = \int_{\mathcal{H}}\left[\scastate{\lambda}\scastate{\mathcal{R}}^w-\left(\int_{\mathcal{H}}\scastate{\lambda}\scastate{\mathcal{R}}^w\xi_i\text{d}V'\right)\scastate{\omega}\xi_j K^{-1}_{ji}\right]\text{d}\scastate{\Phi}\text{d}V',
  \label{fluid_dissipation_energy_increment_non_uniform}
\end{align}
where $i, j=1,2,3$. It follows from \eqref{fluid_dissipation_energy_increment_non_uniform} that the stabilized term for the fluid flow state accounting for the non-uniform fluid potential state can be written as
\begin{equation}
\scastate{\mathcal{Q}}^s = \scastate{\lambda}\scastate{\mathcal{R}}^w-\left(\int_{\mathcal{H}}\scastate{\lambda}\scastate{\mathcal{R}}^w \xi_i\text{d}V'\right)\scastate{\omega}\xi_jK^{-1}_{ji}.
\label{fluid_flow_state_stabilized_term}
\end{equation}
Given $\text{d}\scastate{\Phi}_w$ it can be proved that $\text{d}\mathscr{V}_s$ is always positive for a positive value of $\scastate{\lambda}$. From \eqref{fluid_flow_state_stabilized_term} and the original multiphase correspondence principle \cite{song2020peridynamic} the stabilized fluid flow state can be written as
\begin{equation}
{\scastate{\mathcal{Q}}} = \scastate{\omega}\left(\vec{q}_w-k_w\ddot{\vec{u}}\right)\vec{K}^{-1} \vec{\xi} + \scastate{\lambda}\scastate{\mathcal{R}}^w - \left(\int_{\mathcal{H}}\scastate{\lambda}\scastate{\mathcal{R}}^w\vec{\xi}\text{d}V'\right)\scastate{\omega}\vec{K}^{-1}\vec{\xi}.
\label{fluid_flow_state__stabilization_full_form}
\end{equation}
To simplify \eqref{fluid_flow_state__stabilization_full_form}, it is assumed that
\begin{equation}
\scastate{\lambda} = \dfrac{G K_p}{\omega_0}\scastate{\omega},
\label{ad_hoc_lambda}
\end{equation}
where $G$ is a positive constant on the order of 1 as defined earlier and $K_p$ is a micro-conductivity which will be derived in the next section. With \eqref{ad_hoc_lambda} and \eqref{non_uniform_fluid_potential_gradient}, it can be shown that the third term of \eqref{fluid_flow_state__stabilization_full_form} becomes zero. Thus, equation \eqref{fluid_flow_state__stabilization_full_form} can be expressed as
\begin{equation}
{\scastate{\mathcal{Q}}} = \scastate{\omega}\left[\left(\vec{q}_w-k_w\ddot{\vec{u}}\right)\vec{K}^{-1} \vec{\xi} + \dfrac{G K_p}{\omega_0}\scastate{\mathcal{R}}^w\right].
\label{fluid_flow_state_stabilization_simple_form}
\end{equation}

With \eqref{effective_force_state_stabilization_simple_form} and \eqref{fluid_flow_state_stabilization_simple_form},\eqref{momentum_balance_correspondence} and \eqref{mass_balance_correspondence} with stabilization terms respectively can be rewritten as 
\begin{align}
&\int_\mathcal{H}\left\{\left[\scastate{\omega}(\overline{\vec{\sigma}} - p_w{\vec{1}})J\widetilde{\vec{F}}^{-T}\vec{K}^{-1}\vec{\xi}+  \dfrac{GC}{\omega_0}\vecstate{\mathcal{R}}^s\right] -\left[\scastate{\omega}(\overline{\vec{\sigma}}' - p'_w{\vec{1}})J\widetilde{\vec{F}}^{-T}\vec{K}^{-1}\vec{\xi} + \dfrac{GC}{\omega_0}\vecstate{\mathcal{R}}^{'s}\right]\right\} \text{d}V' \nonumber\\
&+ \rho\vec{g}= \rho\ddot{\vec{u}}, \label{momentum_balance_correspondence_stabilized}\\
&\frac{\phi}{K_w}\dot{p}_w + \int_{\mathcal H} \left(\scastate{\omega}\dot{\vec{u}}_s \vec{K}^{-1}\vec{\xi} - \scastate{\omega}\dot{\vec{u}}'_s \vec{K}^{-1}\vec{\xi}'\right)dV' 
+ \int_{\mathcal{H}} \left\{\scastate{\omega}\left[ (\vec{q}_w - k_w\ddot{\vec{u}})\vec{K}^{-1}\vec{\xi} + \dfrac{G K_p}{\omega_0}\scastate{\mathcal{R}}^w\right] \right.   \nonumber\\
&\left. -  \scastate{\omega} \left[(\vec{q}'_w - k_w\ddot{\vec{u}}')\vec{K}^{-1}\vec{\xi}'+  \dfrac{G K_p}{\omega_0}\scastate{\mathcal{R}}^{' w}\right]\right\}\text{d}V'= 0,
\label{mass_balance_correpondence_stablized}
\end{align}

\subsection{Determination of stabilization parameters}

In this part the parameters in the stabilized terms will be derived based on the energy equivalence between peridynamics and classical theory for the solid skeleton and fluid phase.
For the solid phase, it is postulated that the stored elastic energy in the solid skeleton from peridynamics equals to the elastic energy from the classical poromechanic theory at the same material point $\vec{x}$. 

For simplicity, a microelastic peridynamic model \cite{silling2005meshfree} is adopted to determine the elastic energy in the solid skeleton at material point $\vec{x}$. In the microelastic material model, the effective pairwise force function $\vec{f}$ that material point $\vec{x}'$ imposes on material point $\vec{x}$ is determined from a micropotential $w_s$ as
\begin{equation}
\vec{f}(\vec{\eta},\vec{\xi}) = \dfrac{\vec{\xi} + \vec{\eta}}{|\vec{\xi} + \vec{\eta}|}\dfrac{\partial w}{\partial \eta}(\eta, \xi),
\label{effective_force_density}
\end{equation}
where $\eta = |\vec{\eta}|$ and $\xi=|\vec{\xi}|$.
The micro-potential function is a measure of the elastic strain energy stored in a single bond of the solid skeleton due to its deformation. The total strain energy density at point $\vec{x}$ is expressed as 
\begin{equation}
\mathscr{W} = \frac{1}{2}\int_{\mathcal{H}} w(\eta, \xi) \text{d}V',
\label{peridynamic_strain_energy_generic}
\end{equation}
where the factor of $1/2$ means that each endpoint of a bond between two solid material points owns only half the energy in this bond. Let $f$ be the magnitude of $\vec{f}$ as
\begin{equation}
f(\eta,\xi) = \dfrac{\partial w}{\partial \eta}.
\label{f_magnitude_from_potential}
\end{equation}

We further assume a homogeneous solid skeleton under isotropic extension. It follows
\begin{equation}
\eta = \mathscr{C}_1\xi, 
\end{equation}
where $\mathscr{C}_1$ is a constant for all $\vec{\xi}$. Thus $f$ can be written as
\begin{equation}
f = C\mathscr{C}_1 = C\eta/\xi,
\label{f_magnitude_specific_expression}
\end{equation}
where $C$ is the constant defined previously.
It follows from \eqref{f_magnitude_from_potential} and \eqref{f_magnitude_specific_expression} that 
\begin{align}
w = C\eta^2/(2\xi)= C\lambda^2\xi/2.
\label{micro_potential_expression}
\end{align}
Substituting \eqref{micro_potential_expression} into \eqref{peridynamic_strain_energy_generic} gives
\begin{equation}
\mathscr{W} =  \frac{1}{2}\int_{\mathcal{H}}  \dfrac{1}{2}\left(C\mathscr{C}_1^2\xi\right) \text{d}V'= \frac{1}{2} \int_{0}^{\delta}  \dfrac{1}{2}\left(C\mathscr{C}_1^2\xi\right)  (4\pi\xi^2) \text{d}{\xi} = \frac{\pi C\delta^4}{4}\mathscr{C}_1^2.
\label{peridynamic_elastic_strain_energy}
\end{equation}
The elastic strain energy of the solid skeleton at material point $\vec{x}$ from the classical elastic theory under isotropic deformation is 
\begin{equation}
\widetilde{\mathscr{W}} = \frac{1}{2}\overline{\vec\sigma}:{\vec\epsilon} =  \frac{1}{2}({K^e}{\epsilon_v})(\epsilon_v) = \frac{1}{2}{K^e}(3\mathscr{C}_1)^2 =  \frac{9K^e}{2}\mathscr{C}_1^2,
\label{classical_elastic_strain_energy}
\end{equation}
where $K^e$ is the classical elastic bulk modulus and $\epsilon_v$ is the elastic volumetric strain. Combining \eqref{peridynamic_elastic_strain_energy} and \eqref{classical_elastic_strain_energy} leads to an expression for $C$ under three-dimensional condition as
\begin{equation}
C_{3d} = \frac{18K^e}{\pi\delta^4}.
\label{peridynamic_micromodulus_3d}
\end{equation}

For the fluid phase, the pairwise fluid flow density $f$ at material point $\vec{x}$ in the bond-based peridynamics can be determined from a fluid dissipation micropotential $w^f$ through 
\begin{equation}
f_w(\vec{x},\vec{x}') = \frac{\partial w^f}{\partial \Phi}(\vec{x},\vec{x}'),
\label{peridynamic_pairwise_fluid_flow_density}
\end{equation}
where $\Phi = p'_w - p_w$ is the fluid pressure difference between material points $\vec{x}'$ and $\vec{x}$.
The dissipation micropotential $w^f$ represents the dissipation potential along a bond between two material points  and is  a function of the fluid potential scalar state $\Phi$ of that bond. Then the total dissipation potential at point $\vec{x}$ is a summation over all the micropotentials in the family of this point.
\begin{equation}
\mathscr{V} = \frac{1}{2}\int_{\mathcal{H}} w^f(\vec{x}',\vec{x}) \text{d}V',
\label{peridynamic_fluid_potential_generic}
\end{equation}
where similar to the solid phase the factor of $1/2$ means that each endpoint of a bond between two solid material points owns only half the energy in this bond.
The pairwise fluid flow density at $\vec{x}$ is assumed as 
\begin{equation}
f({\bf x'},{\bf x}) = K_p \frac{\Phi}{|\vec{\xi}|},
\label{peridynamic_pairwise_fluid_flow}
\end{equation}
where $K_p$ is the peridynamic hydraulic micro-conductivity. It follows from \eqref{fluid_flow_density} and \eqref{peridynamic_pairwise_fluid_flow_density} that fluid flow dissipation micro-potential can be written as
\begin{equation}
w^f = \frac{1}{2} K_p \frac{\Phi^2}{|\vec{\xi|}}.
\label{peridynamic_dissipation_micro_potential}
\end{equation}
The peridynamic hydraulic micro-conductivity can be related to the classical hydraulic conductivity by equating the peridynamic fluid dissipation potential to the classical fluid dissipation potential at point $\vec{x}$. For simplicity, we assume a linear pressure field in a body, $p_w =\mathscr{C}_2(\vec{1}\cdot\vec{x})$ for a three-dimensional case. 
Thus, the fluid pressures at material points $\vec{x}$ and $\vec{x}'$ are written as
\begin{equation}
p_w =  \mathscr{C}_2(\vec{1}\cdot\vec{x})\;\;\text{and} \; p'_w = \mathscr{C}_2(\vec{1}\cdot\vec{x}').
\label{linear_pressure_field_at_x_x'}
\end{equation}
It follows that 
\begin{equation}
w^f = K_p  \frac{1}{2}\mathscr{C}^2_2\xi.
\end{equation}
Substituting \eqref{linear_pressure_field_at_x_x'} into \eqref{peridynamic_dissipation_micro_potential} and conducting integration in a sphere with radius $\delta$ centered at $\vec{x}$ lead to
\begin{align}
\mathscr{V} =  \frac{1}{2}\int_0^{\delta} \left( \frac{1}{2}K_p\mathscr{C}^2_2 \xi\right) (4\pi\xi^2)\text{d}\xi
 = \frac{\pi K_p \delta^4}{4}\mathscr{C}^2_2.
 \label{peridynamic_fluid_dissipation_energy}
\end{align}
Assuming a homogeneous body and isotropic fluid flow, the classical fluid dissipation energy through Darcy's law at material point $\vec{x}$ is expressed as
\begin{equation}
\widetilde{\mathscr{V}} = \frac{1}{2}(\vec{\nabla} p_w)k_w\vec{1}({\vec{\nabla} p_w}),
\label{classic_fluid_dissipation_energy_generic}
\end{equation}
where $k_w$ is the hydraulic conductivity of saturated porous media. It follows from \eqref{classic_fluid_dissipation_energy_generic} and \eqref{linear_pressure_field_at_x_x'}  that the classical fluid dissipation energy at material point $\vec{x}$ is
\begin{equation}
\widetilde{\mathscr{V}} = \frac{3}{2}k_w \mathscr{C}^2_2.
\label{classic_fluid_dissipation_energy}
\end{equation}
Combining \eqref{peridynamic_fluid_dissipation_energy} and \eqref{classic_fluid_dissipation_energy} gives the hydraulic micro-conductivity under three-dimensional condition as
\begin{equation}
K_{p, 3d} = \frac{6k_w}{\pi\delta^4}. 
\label{peridynamic_hydraulic_micro_conductivity_3d}
\end{equation}

\section{Numerical implementation}

\subsection{Spatial discretization}

\begin{figure}[h]
\centering
  \includegraphics[width=0.4\textwidth]{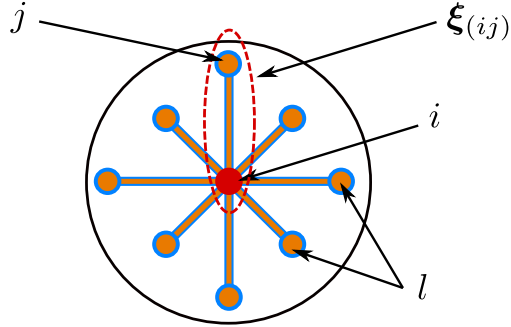}
\caption{Schematic of a Lagrangian-Eulerian meshless spatial discretization scheme of material point $i$ and its neighboring material points (solid points in red and fluid points in blue).}
\label{discretization-scheme}   
\end{figure}

The equations \eqref{momentum_balance_correspondence_stabilized} and \eqref{mass_balance_correpondence_stablized} are spatially discretized by a hybrid Lagrangian-Eulerian meshfree scheme, as shown in Figure \ref{discretization-scheme}.  In this method, a porous continuum material is discretized into a finite number of mixed material points (i.e., mixed solid skeleton and pore water material points). Each material point has two types of degree of freedom, the displacement and the pore water pressure. The uniform grid is used to spatially discretize the problem domain in which all material points have an identical size. 

It is assumed that inertia loading has no impact on the fluid flow \cite{zienkiewicz1999computational} and water is incompressible. Let $\mathscr{N}_i$ be the number of material points in the horizon of a material point $i$. The spatially discretized equations at material point $i$ can be written as
\begin{align}
 \vec{\mathcal{T}}_i + \vec{\mathcal{T}}^s_i+ \rho_i \vec{g} &= \rho_i \ddot{\vec{u}}_i,\label{space_momentum_balance}\\
 \dot{\mathcal{V}}_i + \mathcal{Q}_i + \mathcal{Q}^s_i  &= 0,\label{space_mass_balance}
\end{align} 

where
 \begin{align}
\vec{\mathcal{T}}_i & =\sum_{j = 1}^{\mathscr{N}_i} \left[ \scastate{\omega}J_i\vec{\sigma}_{(i)} \widetilde{\vec{F}}_{(i)}^{-T}\vec{K}_{(i)}^{-1}\vec{\xi}_{(ij)}-  \scastate{\omega}J_j\vec{\sigma}_{(j)} \widetilde{\vec{F}}_{(j)}^{-T}\vec{K}_{(j)}^{-1}\vec{\xi}_{(ji)} \right] V_j , \\
 \vec{\mathcal{T}}^s_i & =\sum_{j = 1}^{\mathscr{N}_i}  \frac{GC}{\omega_0}\left[\vecstate{\mathcal{R}}_{(ij)}^{s} -   \vecstate{\mathcal{R}}_{(ji)}^{s} \right] V_j,\\
\dot{\mathcal{V}}_i  &=  \sum_{j = 1}^{\mathscr{N}_i} \left[\scastate{\omega}(\dot{\vec {u}}_{(j)} - \dot{\vec {u}}_{(i)}){\vec K}_{(i)}^{-1}{\vec \xi}_{(ij)}  \right] V_j, \\
\mathcal{Q}_i  & = \sum_{j = 1}^{\mathscr{N}_i} \left[ \scastate{\omega} {\vec q}_{(i)} {\vec K}_{(i)}^{-1} \vec\xi_{(ij)}  - \scastate{\omega} {\vec q}_{(j)} {\vec K}_{(j)}^{-1} \vec\xi_{(ji)} \right] V_j ,\\
\mathcal{Q}^s_i  & = \sum_{j = 1}^{\mathscr{N}_i} \frac{GK_p}{\omega_0}\left[\scastate{\mathcal{R}}^{w}_{(ij)} -  \scastate{\mathcal{R}}^{'w}_{(ji)} \right] V_j.
b\end{align}

We define the linear assembly operator $\mathscr{A}$ \cite{hughes2012finite}. Let $\mathscr{N}$ be the number of the total material points. The global discretization form of the coupled equations can be written as
\begin{align}
\mathscr{A}_{i=1}^{\mathscr{N}}\left(\vec{\mathcal{T}}_i+ \vec{\mathcal{T}}^s_i+ \rho_i\vec{g}\right)V_i&=\mathscr{A}_{i=1}^{\mathscr{N}} (\rho_i\ddot{\vec{u}}_i)V_i,\\
\mathscr{A}_{i=1}^{\mathscr{N}}\left(\dot{\mathcal{V}}_i+ \mathcal{Q}_i + \mathcal{Q}^s_i\right)V_i  &= \mathscr{A}_{i=1}^{\mathscr{N}}0_i.
\end{align}

\subsection{Time integration and linearization}

We formulate a fully implicit scheme to integrate the coupled system of equations in time\cite{hughes2012finite}. 
At $t_{n+1}$, the residual vector of the the coupled system is defined as
 \begin{align}
&\vec{r}^u_{n+1} =\left.\mathscr{A}_{i=1}^{\mathscr{N}} \left(\rho_i\ddot{\vec{u}}_i- \vec{\mathcal{T}}_i- \vec{\mathcal{T}}^s_i- \rho_i\vec{g}\right)V_i\right|_{n+1} \label{residual_momentum},\\
&\vec{r}^p_{n+1} =\left. \mathscr{A}_{i=1}^{\mathscr{N}}\left(\dot{\mathcal{V}}_i+ \mathcal{Q}_i + \mathcal{Q}^s_i\right)V_i\right|_{n+1} \label{residual_mass}.
\end{align} 

The Newmark method \cite{newmark1959method,zienkiewicz1999computational} is adopted. In the temporal domain a second-order scheme is used to integrate the momentum balance equation and a first-order scheme is applied to integrate the mass balance equation. At time step $n$,  $\vec{u}_n$, $\dot{\vec{u}}_n$, $\ddot{\vec{u}}_n$, $\vec{p}_n$, and $\dot{\vec{p}}_n$ are known. Let $\Delta \ddot{\vec{u}}_{n+1} = \ddot{\vec{u}}_{n+1} - \ddot{\vec{u}}_n$ and $\Delta \dot{\vec{p}}_{n+1} = \dot{\vec{p}}_{n+1} - \dot{\vec{p}}_n$, the acceleration, velocity, displacement and water pressure vectors at $t_{n+1}$ can be written as follows, 
\begin{align}
\ddot{\vec{u}}_{n+1} &= \ddot{\vec{u}}_{n}  + \Delta \ddot{\vec{u}}_{n+1} \label{acceleration_np1},\\
\dot{\vec{u}}_{n+1} &= \dot{\vec{u}}_{n}  + \Delta t  \ddot{\vec{u}}_{n} + \beta_2\Delta t \Delta \ddot{\vec{u}}_{n+1}\label{velocity_np1},  \\
\vec{u}_{n+1} &= \vec{u}_{n} + \Delta t \dot{\vec{u}}_{n} + \frac{(\Delta t)^2}{2}\ddot{\vec{u}}_{n} + \beta_1 \frac{(\Delta t)^2}{2}  \Delta\ddot{\vec{u}}_{n+1}\label{displacement_np1} , \\
\vec{p}_{n+1} &= \vec{p}_n +\Delta t \dot{\vec{p}}_{n} + \beta_3 \Delta t \Delta \dot{\vec{p}}_{n+1}\label{pressure_np1} ,
\end{align}
where $\beta_1, \beta_2, \beta_3\in[0,1]$ are numerical integration parameters. For unconditional stability,
\begin{equation}
\beta_1 \geqslant \beta_2 \geqslant \frac{1}{2}, \quad \mathrm{and} \quad \beta_3 \geqslant \frac{1}{2}.
\end{equation}
Substituting \eqref{acceleration_np1}, \eqref{velocity_np1}, \eqref{displacement_np1}, and \eqref{pressure_np1} into \eqref{residual_momentum} and \eqref{residual_mass}, $\Delta \ddot{\vec{u}}_{n+1} $ and $\Delta \dot{\vec{p}}_{n+1}$, can be solved by Newton's method as follows. Let $k$ be the iteration number.
\begin{equation}
\begin{Bmatrix}
\vec{r}^{u,k+1}\\
\vec{r}^{p,k+1}
\end{Bmatrix}
=\begin{Bmatrix}
\vec{r}^{u,k}\\
\vec{r}^{p,k}
\end{Bmatrix}
+\mathcal{A}^k 
\begin{Bmatrix}
\delta \Delta \ddot{\vec{u}}^{k+1} \\
\delta \Delta \dot{\vec{p}}^{k+1}
\end{Bmatrix}
\approx \begin{Bmatrix}
\vec{0} \\
\vec{0}
\end{Bmatrix}\label{linearization},
\end{equation}
where 
\begin{equation}
\vec{\mathcal{A}}=\left.
\begin{bmatrix}
\dfrac{\partial \vec{r}^u}{\partial \Delta \ddot{\vec{u}}} & \dfrac{\partial \vec{r}^u}{\partial \Delta\dot{\vec{p}}} \\
\vspace{-2ex} \\
\dfrac{\partial \vec{r}^p} {\partial \Delta\ddot{\vec{u}}} & \dfrac{\partial \vec{r}^p}{\partial \Delta\dot{\vec{p}}}
\end{bmatrix}\right|^k_{n+1}
\label{tangent-matrix}.
\end{equation}
By solving \eqref{linearization}, we have
\begin{equation}
\begin{Bmatrix}
\delta \Delta \ddot{\vec{u}}^{k+1} \\
\delta \Delta \dot{\vec{p}}^{k+1}
\end{Bmatrix} = -\vec{\mathcal{A}}^{-1} \begin{Bmatrix}
\vec{r}^{u,k} \\
\vec{r}^{p,k}
\end{Bmatrix}.
\end{equation}
Finally, we have
\begin{equation}
\begin{Bmatrix}
 \Delta \ddot{\vec{u}}^{k+1} \\
 \Delta \dot{\vec{p}}^{k+1}
\end{Bmatrix} = 
\begin{Bmatrix}
\Delta \ddot{\vec{u}}^{k} \\
\Delta \dot{\vec{p}}^{k}
\end{Bmatrix} + 
\begin{Bmatrix}
\delta \Delta \ddot{\vec{u}}^{k+1} \\
\delta \Delta \dot{\vec{p}}^{k+1}
\end{Bmatrix} .
\end{equation}

\subsection{Tangent operator}

Given the relationships in equation \eqref{acceleration_np1} - \eqref{pressure_np1}, by chain rule the individual components of $\vec{\mathcal{A}}$ can be written as follows.
\begin{align}
\bar{\mathbb{K}}_{\vec{uu}} = \dfrac{\partial \vec{r}^u}{\partial \Delta \ddot{\vec{u}}} &= \dfrac{\partial (\vec{\rho}\ddot{\vec{u}})}{\partial \Delta \ddot{\vec{u}}} - \frac{\partial\widehat{\vec{T}} }{\partial \vec{u}}  \frac{\partial \vec{u}}{\partial \Delta \ddot{\vec{u}}} = \vec{\rho}\vec{I}_1 + \dfrac{1}{2}\beta_1\Delta t^2  \left(\frac{\partial \vec{\rho} }{\partial {\vec{u}}}(\ddot{\vec{u}} - \vec{g}) - \frac{\partial\widehat{\vec{T}} }{\partial \vec{u}} \right), \label{momentum-simple-linear-displ}\\
\bar{\mathbb{K}}_{\vec{up}} =\dfrac{\partial \vec{r}^u}{\partial \Delta \dot{\vec{p}}}  &= -\dfrac{\partial \widehat{\vec{T}}}{\partial {\vec{p}}} \frac{\partial \vec{p}}{\partial \Delta \dot{\vec{p}}} = -\beta_3\Delta t \frac{\partial \widehat{\vec{T}}}{\partial \vec{p}}, \label{momentum-simple-linear-pressure} \\
\bar{\mathbb{K}}_{\vec{pu}} = \dfrac{\partial \vec{r}^p}{\partial \Delta \ddot{\vec{u}}} &= \dfrac{\partial \dot{\mathcal{V}}}{\partial \dot{\vec{u}}} \frac{\partial \dot{\vec{u}}}{\partial \Delta \ddot{\vec{u}}} =  \beta_2 \Delta t \frac{\partial \dot{\mathcal{V}}}{\partial \dot{\vec{u}}}\label{mass-balance-linear-displ}, \\
\bar{\mathbb{K}}_{\vec{pp}} = \dfrac{\partial \vec{r}^p}{\partial \Delta \dot{\vec{p}}} &= \left( \frac{\partial \widehat{\vec{\mathcal{Q}}}}{\partial \vec{p} } \right) \frac{\partial {\vec{p}}}{\partial \Delta \dot{\vec{p}}} = \beta_3\Delta t \left( \frac{\partial \widehat{\vec{\mathcal{Q}}}}{\partial \vec{p} } \right). \label{mass-balance-linear-pressure}
\end{align}
where $\bar{\mathbb{K}}_{\vec{up}}$ is the global solid tangent operator, $\bar{\mathbb{K}}_{\vec{up}}$ and $\bar{\mathbb{K}}_{\vec{pu}}$ are the global coupling matrices, $\bar{\mathbb{K}}_{\vec{pp}}$ is the global fluid tangent operator, $\vec{I}_1$ is the second-order identity tensor with the dimension of the number of total material points in the problem domain, and for brevity we define
\begin{align}
\widehat{\vec{T}} &= \vec{\mathcal{T}} + \vec{\mathcal{T}}^s,\\
\widehat{\vec{\mathcal{Q}}} &= \vec{\mathcal{Q}} + \vec{\mathcal{Q}}^s.
\end{align}

Inspired by the standard procedure in the finite element method \cite{hughes2012finite}, the stiffness matrix here will be constructed from the corresponding local stiffness matrices at material points. In this article, the stiffness matrix at a material point will be first computed and then the assembly operator will be utilized to construct the global stiffness matrix. In what follows, the derivation is focused on the stiffness matrices at one material point incorporating all the material points in its horizon.

We linearize the momentum balance equation at material point $i$ by following the chain rule.
The incremental forms of $\vec{\mathcal{T}}_i$ and $\vec{\mathcal{T}}^s_i$ at material point $i$ can be written as
\begin{align}
    \delta \vec{\mathcal{T}}_i &= \sum_{l=1}^{\mathscr{N}_i} \frac{\partial \vec{\mathcal{T}}_i}{\partial \vecstate{Y}_{il}} \frac{\partial \vecstate{Y}_{il}}{\partial \vec{u}_i} \delta\vec{u}_i + \sum_{l=1}^{\mathscr{N}_i}\frac{\partial \vec{\mathcal{T}}_i}{\partial \vecstate{Y}_{il}} \frac{\partial \vecstate{Y}_{il}}{\partial \vec{u}_l} \delta \vec{u}_l +  \frac{\partial \vec{\mathcal{T}}_i}{\partial p_i} \delta p_i + \sum_{l=1}^{\mathscr{N}_i}\frac{\partial \vec{\mathcal{T}}_i}{\partial p_l} \delta p_l \label{total-force-linearization}, \\
    \delta \vec{\mathcal{T}}_i^s &= \sum_{l=1}^{\mathscr{N}_i} \frac{\partial \vec{\mathcal{T}}_i^s}{\partial \vecstate{Y}_{il}} \frac{\partial \vecstate{Y}_{il}}{\partial \vec{u}_i} \delta\vec{u}_i + \sum_{l=1}^{\mathscr{N}_i}\frac{\partial \vec{\mathcal{T}}_i^s}{\partial \vecstate{Y}_{il}} \frac{\partial \vecstate{Y}_{il}}{\partial \vec{u}_l} \delta \vec{u}_l \label{total-force-density-linearization-stab},
\end{align}
where $\vecstate{Y}_{il} = \vec{y}_{l} - \vec{y}_{i}$.

It follows from \eqref{total-force-linearization} and \eqref{total-force-density-linearization-stab} that we have
\begin{align}
 \frac{\partial \vec{\mathcal{T}}_i}{\partial \vec{u}_i}  &= \sum_{l=1}^{\mathscr{N}_i} \frac{\partial \vec{\mathcal{T}}_i}{\partial \vecstate{Y}_{il}} \frac{\partial \vecstate{Y}_{il}}{\partial \vec{u}_i} =  \sum_{l=1}^{\mathscr{N}_i}\left[\sum_{j=1}^{\mathscr{N}_i} \frac{\partial }{\partial \vecstate{Y}_{il}} \left( \vecstate{T}_{ij} -\vecstate{T}_{ji}  \right) V_j \right] \dfrac{\partial \vecstate{Y}_{il}}{\partial \vec{u}_i}, \label{linearize-point-force-point}\\
 \frac{\partial \vec{\mathcal{T}}_i}{\partial \vec{u}_l}  &= \frac{\partial \vec{\mathcal{T}}_i}{\partial \vecstate{Y}_{il}} \frac{\partial \vecstate{Y}_{il}}{\partial \vec{u}_l} = \left[ \sum_{j=1}^{\mathscr{N}_i} \frac{\partial }{\partial \vecstate{Y}_{il}} \left( \vecstate{T}_{ij} - \vecstate{T}_{ji}  \right)V_j \right]\dfrac{\partial \vecstate{Y}_{il}}{\partial \vec{u}_l}\label{linearize-point-force-neighbor},\\
    \frac{\partial \vec{\mathcal{T}}_i^s}{\partial \vec{u}_i} &= \sum_{l=1}^{\mathscr{N}_i} \frac{\partial \vec{\mathcal{T}}_i^s}{\partial \vecstate{Y}_{il}} \frac{\partial \vecstate{Y}_{il}}{\partial \vec{u}_i} =  \sum_{l=1}^{\mathscr{N}_i}\left[\sum_{j=1}^{\mathscr{N}_i}\frac{GC}{\omega_0}  \frac{\partial }{\partial \vecstate{Y}_{il}} \left( \vecstate{\mathcal{R}}^s_{ij} -\vecstate{\mathcal{R}}^s_{ji}  \right)V_j \right] \dfrac{\partial \vecstate{Y}_{il}}{\partial \vec{u}_i}\label{linearize-point-force-point-stab}, \\
    \frac{\partial \vec{\mathcal{T}}_i^s}{\partial \vec{u}_l} &= \frac{\partial \vec{\mathcal{T}}_i^s}{\partial \vecstate{Y}_{il}} \frac{\partial \vecstate{Y}_{il}}{\partial \vec{u}_l} =  \left[\sum_{j=1}^{\mathscr{N}_i}\frac{GC}{\omega_0}  \frac{\partial }{\partial \vecstate{Y}_{il}} \left( \vecstate{\mathcal{R}}^s_{ij} -\vecstate{\mathcal{R}}^s_{ji}  \right)V_j \right] \dfrac{\partial \vecstate{Y}_{il}}{\partial \vec{u}_l}. \label{linearize-point-force-neighbor-stab}
\end{align}
where $\partial \vec{\mathcal{T}}_i/\partial \vecstate{Y}_{il}$ and $\partial \vec{\mathcal{T}}^s_i/\partial \vecstate{Y}_{il}$ can be determined from a peridynamic material model (e.g., ordinary or non-ordinary). For brevity of notations, let us define
\begin{align}
\frac{\partial \widehat{\vec{T}}_i}{\partial \vec{u}_i} &=  \dfrac{\partial \vec{\mathcal{T}}_i}{\partial \vec{u}_i} +  \dfrac{\partial \vec{\mathcal{T}}^{s}_i}{\partial \vec{u}_i}, \label{linearize-point-force-point-combined}\\
\frac{\partial \widehat{\vec{T}}_i}{\partial \vec{u}_l} &= \dfrac{\partial \vec{\mathcal{T}}_i}{\partial \vec{u}_l} +  \dfrac{\partial \vec{\mathcal{T}}^{s}_i}{\partial \vec{u}_l}, \label{linearize-point-force-neighbor-combined}
\end{align} 
For the dynamic loading term, $\rho\ddot{\vec{u}}$, we define,
\begin{align}
\mathbb{M}_i &=  \rho_i\vec{1} + \frac{\beta_1\Delta t^2}{2} \left(\frac{\partial \rho_i}{\partial \vec{u}_i} (\ddot{\vec{u}}_i - \vec{g}) \right),  \label{linearize-inertial-term-point}\\
 \mathbb{M}_l &= \frac{\beta_1\Delta t^2}{2} \left(\frac{\partial \rho_i}{\partial \vec{u}_l} (\ddot{\vec{u}}_i - \vec{g}) \right) \label{linearize-inertial-term-neighbor},
\end{align}
where
\begin{align}
    \frac{\partial \rho_i } {\partial \vec{u}_i} = \rho_s (-\frac{\partial\phi_i}{\partial \vec{u}_i}) + \rho_w \frac{\partial\phi_i}{\partial \vec{u}_i}, \\
    \frac{\partial \rho_i } {\partial \vec{u}_l} = \rho_s (-\frac{\partial\phi_i}{\partial \vec{u}_l}) + \rho_w \frac{\partial\phi_i}{\partial \vec{u}_l}.
\end{align}
Therefore, from \eqref{linearize-point-force-point} - \eqref{linearize-inertial-term-neighbor}, the solid tangent matrix $\bar{\mathbb{K}}^i_{\vec{uu}}$ at material point $i$ can be constructed as
\begin{align}
\bar{\mathbb{K}}^i_{\vec{uu}} &= \begin{bmatrix}
\mathbb{M}_i - \dfrac{\beta_1\Delta t^2}{2} \dfrac{\partial \widehat{\vec{T}}_i}{\partial \vec{u}_i}  & \mathbb{M}_1 - \dfrac{\beta_1\Delta t^2}{2} \dfrac{\partial \widehat{\vec{T}}_i}{\partial \vec{u}_1}   & \cdots & \mathbb{M}_{\mathscr{N}_i} - \dfrac{\beta_1\Delta t^2}{2} \dfrac{\partial \widehat{\vec{T}}_i}{\partial \vec{u}_{\mathscr{N}_i}}
\end{bmatrix}\label{local-stiffness-uu},
\end{align}
where $\bar{\mathbb{K}}^i_{\vec{uu}}$ is a matrix with dimensions ${3 \times 3(\mathscr{N}_i+1)}$. 

Similarly, it follows from \eqref{total-force-linearization} that we have
\begin{align}
\frac{\partial \vec{\mathcal{T}}_i}{\partial p_i} &=\sum_{j=1}^{\mathscr{N}_i} \frac{\partial\vecstate{T}_{ij}  }{\partial p_i} V_j , \label{linearize-point-force-point-pressure}\\
\frac{\partial \vec{\mathcal{T}}_i}{\partial p_j} &= -\frac{\partial\vecstate{T}_{ji} }{\partial p_j} V_j \label{linearize-point-force-neighbor-pressure} .
\end{align}
Thus, the coupling matrix $\bar{\mathbb{K}}^i_{\vec{up}}$ at material point $i$ can be constructed as
\begin{equation}
\bar{\mathbb{K}}^i_{\vec{up}} = -\beta_3 \Delta t \begin{bmatrix}
\dfrac{\partial \vec{\mathcal{T}}_i}{\partial p_i} & \dfrac{\partial \vec{\mathcal{T}}_i}{\partial p_1} & \cdots & \dfrac{\partial \vec{\mathcal{T}}_i}{\partial p_{\mathscr{N}_i}} 
\end{bmatrix}\label{local-stifness-up},
\end{equation}
where $\bar{\mathbb{K}}^i_{\vec{up}}$ is a matrix with dimensions ${3 \times (\mathscr{N}_i+1)}$.

Next, we derive the tangent matrices associated with the mass balance equation at material point $i$.
The linearization of ${\mathcal{Q}}_i$ and ${\mathcal{Q}}^s_i$ at point $i$ can be written as
\begin{align}
    \delta\mathcal{Q}_i = \sum_{l=1}^{\mathscr{N}_i} \frac{\partial \mathcal{Q}_i}{\partial \scastate{\Phi}_{il}} \frac{\partial \scastate{\Phi}_{il}}{\partial {p}_i} \delta{p}_i + \sum_{l=1}^{\mathscr{N}_i}\frac{\partial \mathcal{Q}_i}{\partial \scastate{\Phi}_{il}} \frac{\partial \scastate{\Phi}_{il}}{\partial {p}_l}\delta {p}_l, \label{total-flow-linearization}\\
    \delta\mathcal{Q}_i^s = \sum_{l=1}^{\mathscr{N}_i} \frac{\partial \mathcal{Q}_i^s}{\partial \scastate{\Phi}_{il}} \frac{\partial \scastate{\Phi}_{il}}{\partial {p}_i} \delta{p}_i + \sum_{l=1}^{\mathscr{N}_i}\frac{\partial \mathcal{Q}_i^s}{\partial \scastate{\Phi}_{il}} \frac{\partial \scastate{\Phi}_{il}}{\partial {p}_l}\delta {p}_l\label{total-flow-linearization-stab}.
\end{align}
where $\scastate{\Phi}_{il} = p_l - p_i$.
From \eqref{total-flow-linearization} and \eqref{total-flow-linearization-stab} we obtain
\begin{align}
\frac{\partial \mathcal{Q}_{i}}{\partial p_{i}} &=  \sum_{l=1}^{\mathscr{N}_i} \frac{\partial \mathcal{Q}_i}{\partial \scastate{\Phi}_{il}} \frac{\partial \scastate{\Phi}_{il}}{\partial {p}_i}= \sum_{l=1}^{\mathscr{N}_i}\left[\sum_{j=1}^{\mathscr{N}_i} \frac{\partial }{\partial \scastate{\Phi}_{il}} \left( \scastate{\mathcal{Q}}_{ij} - \scastate{\mathcal{Q}}_{ji}  \right)V_j \right]\frac{\partial \scastate{\Phi}_{il}}{\partial p_i}, \label{linearize-point-flux-point}\\
\frac{\partial \mathcal{Q}_{i}}{\partial p_{l}} &=\frac{\partial \mathcal{Q}_i}{\partial \scastate{\Phi}_{il}} \frac{\partial \scastate{\Phi}_{il}}{\partial {p}_l}=  \left[\sum_{j=1}^{\mathscr{N}_i} \frac{\partial }{\partial \scastate{\Phi}_{il}} \left( \scastate{\mathcal{Q}}_{ij} - \scastate{\mathcal{Q}}_{ji} \right) V_j \right]\frac{\partial \scastate{\Phi}_{il}}{\partial p_l}\label{linearize-point-flux-neighbor},\\
\frac{\partial \mathcal{Q}_{i}^s}{\partial p_{i}} &= \sum_{l=1}^{\mathscr{N}_i} \frac{\partial \mathcal{Q}^s_{i}}{\partial \scastate{\Phi}_{il}} \frac{\partial \scastate{\Phi}_{il}} {\partial p_{i}} = \sum_{l=1}^{\mathscr{N}_i}\left[\sum_{j=1}^{\mathscr{N}_i} \frac{GK_p}{\omega_0} \frac{\partial }{\partial \scastate{\Phi}_{il}} \left( \scastate{\mathcal{R}}^w_{ij} - \scastate{\mathcal{R}}^w_{ji} \right) V_j \right]\frac{\partial \scastate{\Phi}_{il}}{\partial p_i}, \label{linearize-point-flux-point-stab}\\
\frac{\partial \mathcal{Q}^s_{i}}{\partial p_{l}} &= \frac{\partial \mathcal{Q}^s_{i}}{\partial \scastate{\Phi}_{il}} \frac{\partial \scastate{\Phi}_{il}} {\partial p_{l}} = \left[\sum_{j=1}^{\mathscr{N}_i}\frac{GK_p}{\omega_0}  \frac{\partial }{\partial \scastate{\Phi}_{il}} \left( \scastate{\mathcal{R}}^w_{ij} - \scastate{\mathcal{R}}^w_{ji} \right) V_j \right]\frac{\partial \scastate{\Phi}_{il}}{\partial p_l}\label{linearize-point-flux-neighbor-stab}.
\end{align}
For brevity of notations, let us define 
\begin{align}
\dfrac{\partial \widehat{\mathcal{Q}}_i}{\partial p_i} = \dfrac{\partial \mathcal{Q}_i}{\partial p_i} + \dfrac{\partial \mathcal{Q}^{s}_i}{\partial p_i} \label{linearize-point-flux-point-combined},\\
\dfrac{\partial \widehat{\mathcal{Q}}_i}{\partial p_{l}} = \dfrac{\partial \mathcal{Q}_i}{\partial p_{l}} + \dfrac{\partial \mathcal{Q}^{s}_i}{\partial p_{l}}\label{linearize-point-flux-neighbor-combined}.
\end{align}
It follows from \eqref{linearize-point-flux-point-combined} and \eqref{linearize-point-flux-neighbor-combined} that the flow tangent matrix $\bar{\mathbb{K}}^i_{\vec{pp}}$  at point $i$ can be written as
\begin{equation}
\bar{\mathbb{K}}^i_{\vec{pp}} = \beta_3\Delta t \begin{bmatrix}
\dfrac{\partial \widehat{\mathcal{Q}}_{i}}{\partial p_{i}} & \dfrac{\partial \widehat{\mathcal{Q}}_{i}}{\partial p_{1}} & \cdots & \dfrac{\partial \widehat{\mathcal{Q}}_{i}}{\partial p_{\mathscr{N}_i}}
\end{bmatrix} \label{local-stifness-pp},
\end{equation}
where $\bar{\mathbb{K}}^i_{\vec{pp}}$ is a row vector with the dimension ${(\mathscr{N}_i+1)}$.
For the rate of the solid volume change, it can be readily derived by chain rule that
\begin{align}
    \frac{\partial \dot{\mathcal{V}}_i}{\partial \dot{\vec{u}}_{i}} &= -\sum_{j = 1}^{\mathscr{N}_i} \left[\scastate{\omega}\vec{1}{\vec K}_{(i)}^{-1}{\vec \xi}_{(ij)} \right] V_j, \label{linearize-volume-rate-point} \\
    \frac{\partial \dot{\mathcal{V}}_i}{\partial \dot{\vec{u}}_{j}} &= \left[\scastate{\omega}\vec{1}{\vec K}_{(i)}^{-1}{\vec \xi}_{(ij)} \right] V_j. \label{linearize-volume-rate-neighbor}
\end{align}
Thus, the coupling tangent matrix $\bar{\mathbb{K}}^i_{\vec{pu}}$ at material point $i$ can be written as
\begin{equation}
\bar{\mathbb{K}}^i_{\vec{pu}}= \beta_2 \Delta t
\begin{bmatrix}
\dfrac{\partial \dot{\mathcal{V}}_i}{\partial \dot{\vec{u}}_{i}} & \dfrac{\partial \dot{\mathcal{V}}_i}{\partial \dot{\vec{u}}_{1}} & \cdots & \dfrac{\partial \dot{\mathcal{V}}_i}{\partial  \dot{\vec{u}}_{\mathscr{N}_i}}
\end{bmatrix} \label{local-stifness-pu} ,
\end{equation}
where $\bar{\mathbb{K}}^i_{\vec{pu}}$ is a vector with the dimensions of $1\times{3(\mathscr{N}_i+1)}$.

Finally, the global tangent matrix can be constructed from \eqref{local-stifness-uu}, \eqref{local-stifness-pp}, \eqref{local-stifness-up}, and \eqref{local-stifness-pu} for all material points.
\begin{align}
\bar{\mathbb{K}}_{\vec{uu}} &= \mathlarger{\mathscr{A}}^{\mathscr{P}}_{i=1} \bar{\mathbb{K}}^i_{\vec{uu}} ,\label{momentum-disp-tangent}\\
\bar{\mathbb{K}}_{\vec{up}} &= \mathlarger{\mathscr{A}}^{\mathscr{P}}_{i=1} \bar{\mathbb{K}}^i_{\vec{up}}  \label{momentum-press-tangent},\\
\bar{\mathbb{K}}_{\vec{pu}} &=  \mathlarger{\mathscr{A}}^{\mathscr{P}}_{i=1} \bar{\mathbb{K}}^i_{\vec{pu}}, 
\label{mass-disp-tangent}\\
\bar{\mathbb{K}}_{\vec{pp}}  & = \mathlarger{\mathscr{A}}^{\mathscr{P}}_{i=1} \bar{\mathbb{K}}^i_{\vec{pp}}  \label{mass-press-tangent},
\end{align}
where $\mathlarger{\mathscr{A}}$ is the assembly operator that constructs the global stiffness matrix from the tangent matrices at individual point. In the computer code, each material point and material points in its horizon are assigned a unique global identification (ID) number that are used by $\mathlarger{\mathscr{A}}$ to construct the global tangent matrix. Indeed the assembly procedure adopted here is exactly the same as the global stiffness assembly operator in the finite element method \cite{hughes2012finite}. Algorithm 1 summarize the general procedure for the computation of local tangent matrices at the material point level and the assembly of the global tangent matrix. Parallel computing through Open MPI \cite{gabriel2004open} is exploited to increase the computational efficiency of the stabilized coupled periporomechanics model.

\begin{algorithm}
\setstretch{1.0}
\caption{Compute local tangent matrices and assembly the global tangent matrix}
\begin{algorithmic}[1]
\Procedure{assemble tangent matrix}{}
	\State Allocate and initialize the global tangent matrix.
	\For {node $i \in \mathscr{P}$}
		\For {neighbor $j \in \mathscr{N}_i$}
			\State Compute ${\partial \vec{\mathcal{T}}_i}/{\partial \vec{u}_i}$ and ${\partial \vec{\mathcal{T}}^{s}_i}/{\partial \vec{u}_i}$ using \eqref{linearize-point-force-point} and \eqref{linearize-point-force-point-stab}.
            \State Compute ${\partial \mathcal{Q}_{i}}/{\partial p_{i}}$ and ${\partial \mathcal{Q}^{s}_i}/{\partial p_{i}}$ using \eqref{linearize-point-flux-point} and \eqref{linearize-point-flux-point-stab}.
            
            \State Compute ${\partial\vec{\mathcal{T}}_{i}}/{\partial {p}_i}$ using \eqref{linearize-point-force-point-pressure}.

            \State Compute ${\partial \dot{\mathcal{V}}_i}/{\partial \dot{\vec{u}}_{i}}$ using \eqref{linearize-volume-rate-point}.
            
            \State Compute ${\partial \dot{\mathcal{V}}_i}/{\partial \dot{\vec{u}}_{j}}$ using \eqref{linearize-volume-rate-neighbor}.
            
            \State Compute $\mathbb{M}_j$ using \eqref{linearize-inertial-term-neighbor}.

			\For {neighbor $l \in \mathscr{N}_i$}
			    \State Compute ${\partial \vec{\mathcal{T}}_i}/{\partial \vec{u}_l}$ and ${\partial \vec{\mathcal{T}}^{s}_i}/{\partial \vec{u}_l}$ using \eqref{linearize-point-force-neighbor} and \eqref{linearize-point-force-neighbor-stab}.
				\If{$j == l$}
				    \State\quad Compute ${\partial\vec{\mathcal{T}}_{i}}/{\partial {p}_l}$ using \eqref{linearize-point-force-neighbor-pressure}.
				\EndIf
				\State Compute ${\partial \mathcal{Q}_{i}}/{\partial p_{l}}$ and ${\partial \mathcal{Q}^{s}_i}/{\partial p_{l}}$ using \eqref{linearize-point-flux-neighbor} and \eqref{linearize-point-flux-neighbor-stab}.
			\EndFor			
		\EndFor
		\State Compute $\mathbb{M}_i$ using \eqref{linearize-inertial-term-point}.
		\State Construct $\bar{\mathbb{K}}^i_{\vec{uu}}$, $\bar{\mathbb{K}}^i_{\vec{up}}$, $\bar{\mathbb{K}}^i_{\vec{pu}}$ and $\bar{\mathbb{K}}^i_{\vec{pp}}$.

	\EndFor
	\State Assemble the global tangent matrix using \eqref{momentum-disp-tangent} - \eqref{mass-press-tangent}.
\EndProcedure

\end{algorithmic}
\end{algorithm}

\subsection{Linearization of correspondence material models}

We present the material-point level linearization of correspondence solid and fluid models implemented in this study. Now, the differentiation of the total force density at a material point $i$ with respect to a deformation state $\vecstate{Y}_{il}$ reads,
\begin{align}
\frac{\partial \vec{\mathcal{T}}_{i}}{\partial \vecstate{Y}_{il}} &= \sum_{j=1}^{\mathscr{N}_i} \frac{\partial }{\partial \vecstate{Y}_{il}} \left( \vecstate{T}_{ij} -\vecstate{T}_{ji}  \right) V_j \nonumber\\
&= \sum_{j=1}^{\mathscr{N}_i} \left[ \scastate{\omega}\left( \left\{ \frac{\partial\vec{P}_{(i)}}{\partial \widetilde{\vec{F}}_{(i)}} \frac{\partial \widetilde{\vec{F}}_{(i)}}{\partial \vecstate{Y}_{il}} \right\} (\vec{K}_{(i)})^{-1}\vec{\xi}_{ij}\right) - \scastate{\omega} \left(\left\{ \frac{\partial\vec{P}_{(j)}}{\partial \widetilde{\vec{F}}_{(j)}}\frac{\partial \widetilde{\vec{F}}_{(j)}}{\partial \vecstate{Y}_{il}} \right\} (\vec{K}_{(j)})^{-1}\vec{\xi}_{ji} \right) \right]V_j, 
\end{align}
where $\frac{\partial\vec{P}_{(i)}}{\partial \widetilde{\vec{F}}_{(i)}}$ can be determined through the return mapping algorithm following the lines in computational plasticity (e.g., \cite{simo1998computational, song2014mathematical, borja2013shear}, and others)
\begin{align}
\frac{\partial \widetilde{\vec{F}}_{(i)}}{\partial\vecstate{Y}_{il} } &= \scastate{\omega} \left[ \frac{\partial \vecstate{Y}_{il}}{\partial \vecstate{Y}_{il}} \otimes \vec{\xi}_{il} V_l \right]  (\vec{K}_{(i)})^{-1} =  \scastate{\omega} \left[ \vec{1} \otimes \vec{\xi}_{il} V_l \right]  (\vec{K}_{(i)})^{-1} , \\
\frac{\partial \widetilde{\vec{F}}_{(j)}}{\partial\vecstate{Y}_{il} } &= \scastate{\omega} \left[ \frac{\partial \vecstate{Y}_{il}}{\partial \vecstate{Y}_{li}} \otimes \vec{\xi}_{li} V_i \right] (\vec{K}_{(j)})^{-1}= \scastate{\omega} \left[ -\vec{1}\otimes \vec{\xi}_{li} V_i \right] (\vec{K}_{(j)})^{-1}.
\end{align}
The linearization of the stabilization force with respect to a deformation state $\vecstate{Y}_{il}$ reads,
\begin{equation}
\frac{\partial \vec{\mathcal{T}}^{s}_i}{\partial \vecstate{Y}_{il}} =  \sum_{j=1}^{\mathscr{N}_i} \left[\scastate{\omega}\frac{GC}{\omega_0} \left\{\frac{\partial \vecstate{Y}_{ij}}{\partial \vecstate{Y}_{il}}  - \left(\frac{\partial \widetilde{\vec{F}}_{(i)}}{\partial \vecstate{Y}_{il}} \vec{\xi}_{ij}\right)\right\} - \scastate{\omega}\frac{GC}{\omega_0} \left\{\frac{\partial \vecstate{Y}_{ji}}{\partial \vecstate{Y}_{il}}  - \left(\frac{\partial \widetilde{\vec{F}}_{(j)}}{\partial \vecstate{Y}_{il}} \vec{\xi}_{ji}\right)\right\} \right] V_j,
\end{equation}
where $\frac{\partial \vecstate{Y}_{ij}}{\partial \vecstate{Y}_{il}}$ and $\frac{\partial \vecstate{Y}_{ji}}{\partial \vecstate{Y}_{il}}$ is nonzero only if $j=l$.
Next, the momentum balance linearized with respect to the fluid pressure reads,
\begin{align}
    \frac{\partial \vec{\mathcal{T}}_i}{\partial {p}_i} &= \sum_{j=1}^{\mathscr{N}_i} \frac{\partial\vecstate{T}_{ij}}{\partial p_i} = \sum_{j=1}^{\mathscr{N}_i} \scastate{\omega} (\vec{K}_{(i)})^{-1} \vec{\xi}_{ij} V_j, \\
    \frac{\partial \vec{\mathcal{T}}_i}{\partial {p}_j} &= -\frac{\partial\vecstate{T}_{ji}}{\partial p_j} = -\scastate{\omega} (\vec{K}_{(j)})^{-1} \vec{\xi}_{ji} V_j.    
\end{align}

Similarly, for the mass balance, the derivative of fluid flow density with respect to a pressure potential state reads,
\begin{align}
 \frac{\partial {\mathcal{Q}}_{i}}{\partial \scastate{\Phi}_{il}} &= \sum_{j=1}^{\mathscr{N}_i} \frac{\partial }{\partial \scastate{\Phi}_{il}} \left(\scastate{\mathcal{Q}}_{ij} - \scastate{\mathcal{Q}}_{ji} \right)V_j \nonumber \\
 &= \sum_{j=1}^{\mathscr{N}_i} \left[\scastate{\omega} \left\{ \frac{\partial \vec{q}_{i,w}}{\partial \widetilde{\grad{\Phi}}_{i}} \frac{\partial \widetilde{\grad{\Phi}}_{i}}{\partial \scastate{\Phi}_{il} } \right\} (\vec{K}_{(i)})^{-1}\vec{\xi}_{ij} - \scastate{\omega} \left\{ \frac{\partial \vec{q}_{j,w}}{\partial \widetilde{\grad{\Phi}}_{j}} \frac{\partial \widetilde{\grad{\Phi}}_{j}}{\partial \scastate{\Phi}_{il} } \right\} (\vec{K}_{(j)})^{-1}\vec{\xi}_{ji} \right] V_j ,
\end{align}
where
\begin{align}
 \frac{\partial \widetilde{\grad{\Phi}}_i}{\partial \scastate{\Phi}_{il}} &=  \left[ \scastate{\omega} \frac{\partial \scastate{\Phi}_{il}}{\partial \scastate{\Phi}_{il}}  \vec{\xi}_{il} V_j \right] (\vec{K}_{(i)})^{-1} = \scastate{\omega} \vec{\xi}_{il}  (\vec{K}_{(i)})^{-1} V_j , \\
 \frac{\partial \widetilde{\grad{\Phi}}_j}{\partial \scastate{\Phi}_{il}} &=  \left[ \scastate{\omega} \frac{\partial \scastate{\Phi}_{li}}{\partial \scastate{\Phi}_{il}}  \vec{\xi}_{li} V_i \right] (\vec{K}_{(j)})^{-1} = -\scastate{\omega} \vec{\xi}_{li}   (\vec{K}_{(j)})^{-1} V_i.
\end{align}

The linearization of the stabilization flow density with respect to a deformation state $\scastate{\Phi}_{il}$ reads
\begin{equation}
 \frac{\partial {\mathcal{Q}}^{s}_i}{\partial \scastate{\Phi}_{il}} =\sum_{j=1}^{\mathscr{N}_i} \left[ \scastate{\omega} \frac{GK_p}{\omega_0} \left\{ \frac{\partial \scastate{\Phi}_{ij}}{\partial \scastate{\Phi}_{il}} -  \left(\frac{\partial \widetilde{\grad{\Phi}}_{i}}{\partial \scastate{\Phi}_{il}} \vec{\xi}_{ij}\right) \right\} - \scastate{\omega} \frac{GK_p}{\omega_0} \left\{ \frac{\partial \scastate{\Phi}_{ji}}{\partial \scastate{\Phi}_{il}} -  \left(\frac{\partial \widetilde{\grad{\Phi}}_{j}}{\partial \scastate{\Phi}_{il}} \vec{\xi}_{ji}\right) \right\} \right] V_j,
\end{equation}
where $\frac{\partial \scastate{\Phi}_{ji}}{\partial \scastate{\Phi}_{il}}$ and $\frac{\partial \scastate{\Phi}_{ji}}{\partial \scastate{\Phi}_{il}}$ are nonzero only if $j=l$.

\section{Numerical examples}

For the numerical examples in this section we adopt the time integration parameters in the literature that satisfy the criteria for unconditional stability and high frequency numerical dissipation \cite{zienkiewicz1999computational, popescu1993centrifuge}. For all three examples in what follows we have $\beta_1 = 0.605 $ and $\beta_2=\beta_3=0.6$.  All boundary conditions are imposed through a fictitious boundary layer (see \cite{Silling2000Reformulation, Silling2007Peridynamic, madenci2014peridynamic}).

\subsection{Correspondence material models}
For the material models, we implement two classical material models for the solid skeleton and the classical Darcy's law for fluid flow through the extended correspondence principle with stabilization. Following the lines in continuum mechanics \cite{simo1998computational}, the total strain $\vec{\varepsilon}$ can be determined from the nonlocal deformation gradient at a material point.
In examples 1 and 2, the solid skeleton is modeled by an isotropic elastic constitutive relationship, which reads
\begin{equation}
\overline{\sigma}_{ij} = \mathcal{C}^e_{ijkl}\varepsilon^e_{kl},
\label{elastic-model}
\end{equation}
where $i,j,k,l=1,2,3$, $\mathcal{C}^e_{ijkl}$ is the isotropic elastic tensor, and $\varepsilon^e_{ij}$ is the elastic strain tensor. The isotropic elastic tensor is written as
\begin{equation}
\mathcal{C}^e_{ijkl} = K^e\delta_{ij}\delta_{kl} + 2\mu^e (\mathcal{I}_{ijkl} - 1/3\delta_{ij}\delta_{kl}),
\end{equation}
where $\mu^e$ is the shear modulus and $\mathcal{I}_{ijkl}$ is the rank-four identity tensor. 

In example 3, the critical-state based elastoplastic model for saturated soils \cite{schofield1968critical, wood1990soil} is adopted and numerically implemented through the celebrated return-mapping algorithm in computational plasticity \cite{simo1998computational, borja2013plasticity}. In this material model, the total strain tensor is additively decomposed into the elastic strain tensor $\vec{\varepsilon}^e$ and the plastic strain tensor $\vec{\varepsilon}^p$ as
\begin{equation}
\varepsilon_{ij} = {\varepsilon}^e_{ij} + {\varepsilon}^p_{ij},
\label{strain_decomposition}
\end{equation}
Given the elastic strain tensor, the effective stress can be computed by equation \eqref{elastic-model}. The yield function $f$ is written as
\begin{equation}
f(\overline{p},q, p_c)=(\overline{p}-p_c)\overline{p} + (q/M)^2 \leqslant 0,
\end{equation}
where $\overline{p}$ is the effective mean stress, $q$ is the equivalent shear stress, $M$ is the slope of the critical state line, and ${p}_c$ is the preconsolidation pressure. As a hardening law, $p_c$ evolves with plastic volumetric strain.
\begin{equation}
 \dot{p}_c=\dfrac{-p_c}{{\lambda}-{\kappa}}(\dot{\varepsilon}^p_{11} + \dot{\varepsilon}^p_{22} + \dot{\varepsilon}^p_{33} ),
 \end{equation} 
where $\lambda$ and ${\kappa}$ are the swelling index and the compression index, respectively. The plastic strain is determined below by assuming the associative flow rule.
\begin{equation}
\dot{{\varepsilon}}^p_{ij} = \dot{\gamma}\dfrac{\partial f}{\partial\overline{{\sigma}}_{ij}},
\end{equation}
where $\dot{\gamma}$ is the non-negative plastic multiplier. It is noted that advanced constitutive models for geomaterials can be applied in the formulated peridynamics poromechanics in this article by the recently proposed multiphase correspondence principle (see \cite{song2020peridynamic}).

\subsection{One-dimensional dynamic consolidation problem}
\label{Example_1}

This example concerns the one-dimensional dynamic consolidation of a saturated soil specimen under instantaneous and sinusoidal loading conditions, respectively. Numerical results from the coupled peridynamics (PD) model are compared with the analytical and finite element (FE) solutions in the literature \cite{schanz2000transient, de1993one}. 

\begin{figure}[h!]
\centering
  \includegraphics[width=0.3\textwidth]{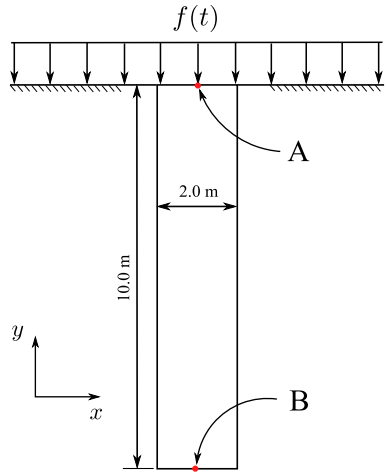}
\caption{Sketch of the soil column and the loading protocol.}
\label{example1-setup}   
\end{figure}

Figure \ref{example1-setup} depicts the geometry of the soil column. The load is imposed on the top boundary. The soil column is restricted to deform vertically while the bottom is fixed. For the fluid phase, the top boundary is drainage and all other boundaries are impervious. The initial effective stress and water pressure are not taken into account as assumed in \cite{schanz2000transient, de1993one}. The problem domain is discretized into 25,000 mixed material points. The center-to-center distance of two neighboring material points is $\Delta x$ = 0.04 m.

\subsubsection{Instantaneous load}

\begin{figure}[h]
\centering
  \includegraphics[width=0.4\textwidth]{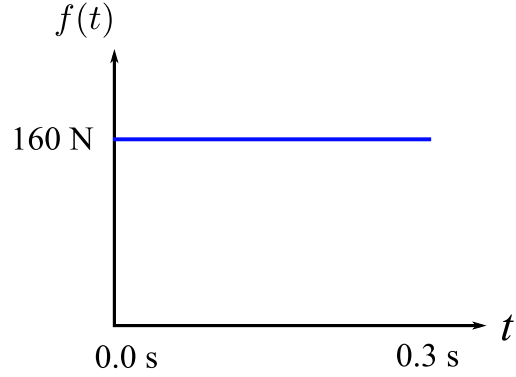}
\caption{Graphical depiction of the instantaneous loading profile.}
\label{example1-instant-load}   
\end{figure}

Figure \ref{example1-instant-load} plots the instantaneous load imposed on the top of the specimen. The simulation time $t=0.3$ s and the time increment $\Delta t$ = 1$\times 10^{-4}$ s \cite{schanz2000transient}. The solid skeleton is modeled by a correspondence elastic model. The fluid flow is modeled by the correspondence Darcy's law. The material parameters adopted from \cite{schanz2000transient} are: bulk modulus $K = 2.1 \times 10^5$ kPa, shear modulus $\mu_s = 9.8 \times 10^4$ kPa, $\rho_s = 1884$ kg/m$^3$, initial porosity $\phi_0 = 0.48$, $K_w$ = $2.2 \times 10^6$ kPa, $\rho_w = 1000$ kg/m$^3$, $k_w = 3.55 \times 10^{-5}$ m/s. The horizon $\delta = 2.05 \Delta x$. The simulation time $t=0.3$ s and the time increment $\Delta t = 1\times 10^{-4}$ s.

\begin{figure}[h!]
\centering
  \includegraphics[width=0.4\textwidth]{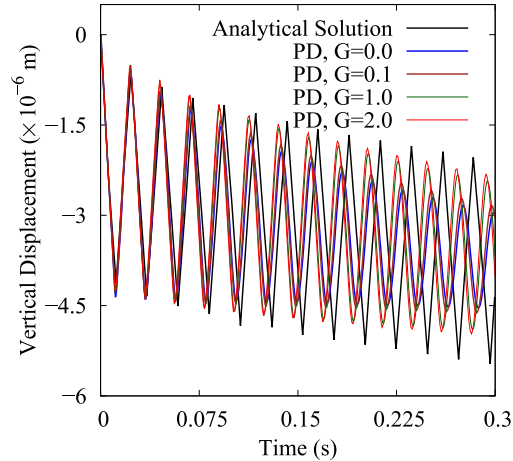}
\caption{Comparison of the PD solution and analytical solution of the vertical displacement at point A.}
\label{1d-bench-results-displacement}   
\end{figure}

\begin{figure}[h!]
\centering
  \includegraphics[width=0.4\textwidth]{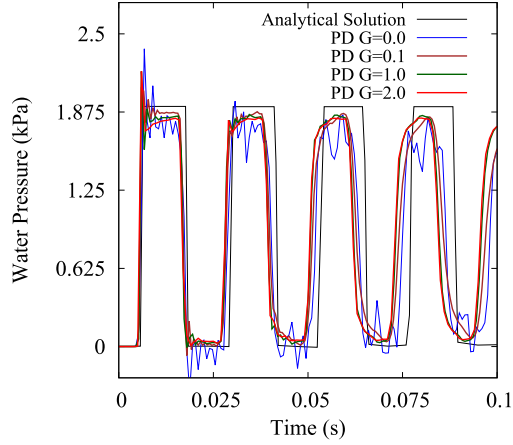}
\caption{Comparison of the PD solution and analytical solution of the water pressure at B.}
\label{1d-bench-results-pressure}   
\end{figure}

Figure \ref{1d-bench-results-displacement} plots the coupled PD solutions with different values of $G$ parameter and analytical solution of the vertical displacement at A shown in Figure \ref{example1-instant-load}. Figure \ref{1d-bench-results-pressure} compares the PD solutions with different values of $G$ parameter and analytical solution of the water pressure at B shown in Figure \ref{example1-instant-load}.

It is shown from Figure \ref{1d-bench-results-displacement} that the value of $G$ can have a significant effect on the predicted amplitude of the oscillations of the vertical displacement at point A. For $G$ = 0, the PD results of vertical displacement are in good agreement with the analytical solution early in the simulation, with the triangular waveform of the analytical solution being largely preserved. These oscillations in amplitude diminish over time due to fluid viscous dissipation, and will eventually go to zero. However, the oscillations in the PD solution appear to dissipate more rapidly than the analytical solution. With $G>0$, the PD solution approaches the analytical one, with $G$ = 1.0 giving the best results. It is also apparent that increasing the value of $G>1$, appears to have a negligible influence on the PD solutions, with practically no difference in the results obtained for $G$ = 1 and $G$ = 2.

The water pressure response at B as shown in Figure \ref{1d-bench-results-pressure} takes the form of periodic square waves of a constant amplitude. For $G$ = 0, the PD solution has noticeable high-frequency oscillations at the peak and valley of the square waves. It is apparent that $G = 0.1$ almost eliminates high-frequency oscillations in the plot of water pressure. Increasing the value of $G$ to 1 or 2 does not improve the accuracy of the PD solution with respect to the analytical one. Indeed, a value of $G$ = 2 may decrease the period of the PD wave shown in the analytical wave profile. 

\subsubsection{Sinusoidal load}

\begin{figure}[h!]
\centering
  \includegraphics[width=0.4\textwidth]{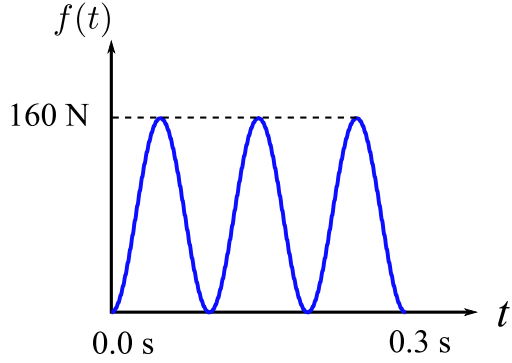}
\caption{Graphical depiction of the sinusoidal loading profile.}
\label{example1-sine-load}   
\end{figure}

The sinusoidal loading profile as shown in Figure \ref{example1-sine-load} is expressed as
\begin{equation}
f (t) = 160 ~ [1-\cos(\omega t)] \label{Harmonic-load},
\end{equation}
where the angular frequency $\omega$ = 20$\pi$ radians/s. The total simulation time $t$ = 0.5 s and the time increment $\Delta t$ = 1$\times 10^{-3}$ s. The solid skeleton is modeled by a correspondence elastic model and the fluid flow is modeled by the correspondence Darcy's law. The input material parameters adopted from \cite{markert2010comparison} are: $K = 1.22 \times 10^4$ kPa, $\mu_s = 5.62 \times 10^3$ kPa, $\rho_s = 2000$ kg/m$^3$, $\phi_0 = 0.33$, $K_w = 2.2 \times 10^6$ kPa, $\rho_w = 1000$ kg/m$^3$, and $k_w = 1.0 \times 10^{-2}$ m/s. The horizon $\delta = 2.05 \Delta x$. The initial and boundary conditions are identical to the instantaneous loading scenario.

\begin{figure}[h!]
\centering
  \includegraphics[width=0.4\textwidth]{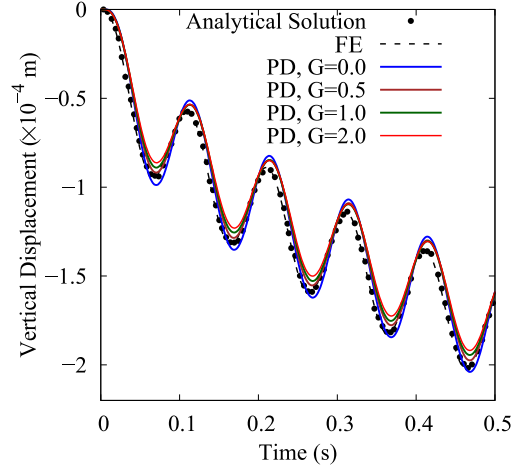}
\caption{Comparison of the PD solution and analytical and FE solutions of the vertical displacement at point A.}
\label{bench1-harmonic-results-displacement}   
\end{figure}

\begin{figure}[h!]
\centering
  \includegraphics[width=0.4\textwidth]{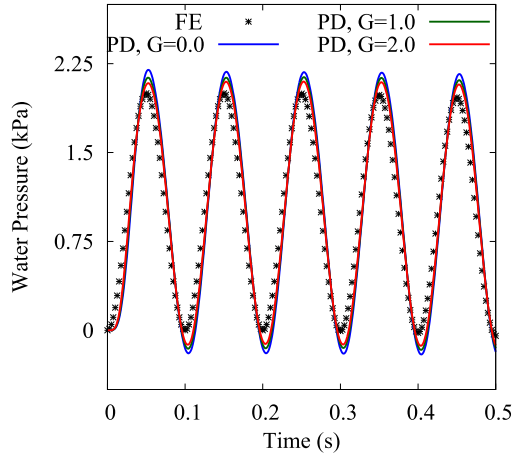}
\caption{Comparison of the PD solution and analytical and FE solutions of the water pressure at point B.}
\label{bench1-harmonic-results-pressure}   
\end{figure}

Figure \ref{bench1-harmonic-results-displacement} plots the vertical displacement at the point A in 
Figure \ref{example1-setup} over time from the PD solution and the analytical and FE solutions. Figure \ref{bench1-harmonic-results-pressure} plots the water pressure at the point B in 
Figure \ref{example1-setup} over time from the PD solution and the FE solution.
The results in Figure \ref{bench1-harmonic-results-displacement} show that the PD solution generally matches the analytical and FE solutions of the vertical displacement at the point A. It is found that the value of $G$ slightly affects the PD solution in this case. With $G$ = 0, the PD solution of the displacement is close to the analytical solution at the peaks of the displacement profile. However, the PD solution diverges from the analtyical solution at the valleys of the displacement profile. With $G > 0$ the PD solution generated a uniform decrease in the amplitude of the oscillations of the displacement profile. Figure \ref{bench1-harmonic-results-pressure} shows that the PD solution with of the water pressure with $G$ = 1.0 is close to the FE solution. There is no noticeable change in the PD solution with $G$ = 2.0.

The influence of the stabilization parameter on the solid skeleton and pore fluid discussed in this example is consistent with the results obtained from the quasi-static analysis of the solid material \cite{hashim2020implicit}. It was found that for simple elastic problems $G$ values as low as 0.01 could provide maximum reduction in relative errors with the analytical solution. Increasing the value of G may gradually lead to an increase of the relative error over the non-stabilized case. Larger values of G may be required when modeling fracture propagation or finite deformation. In the dynamic analysis of solids in \cite{silling2017stability}, it was suggested the value of $G$ should be in the order of 1 for high strain rates and large deformation.  

\subsection{Wave propagation in saturated soil}
\label{Example_2}

This example deals with a two-dimensional saturated soil specimen under dynamic loading through a strip footing. The numerical results from PD are compared with the FE solutions in the literature \cite{markert2010comparison}.

\begin{figure}[h!]
\centering
  \includegraphics[width=0.5\textwidth]{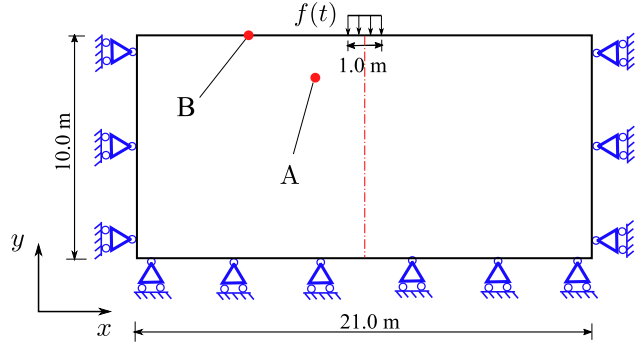}
\caption{Sketch of the problem domain, boundary condition, and loading protocol. }
\label{2d-bench}   
\end{figure}

Figure \ref{2d-bench} depicts the problem geometry and the boundary conditions. For the fluid phase, the fluid pressure at the top boundary is set to zero. The remaining boundaries are impermeable. The solid skeleton is modeled by a correspondence elastic model and the fluid flow is modeled by the correspondence Darcy's law. The problem domain is discretized into 100,000 uniform mixed material points with $\Delta x$ = 0.05 m. $\delta = 2.05 \Delta x$ m. The input material parameters adopted from \cite{pastor2000fractional, breuer1999quasi, markert2010comparison} are: $K = 1.22 \times 10^4$ kPa, $\mu_s = 5.62 \times 10^3$ kPa, $\rho_s = 2000$ kg/m$^3$, $\phi = 0.33$, $K_w$ = $2.2 \times 10^6$ kPa, $\rho_w = 1000$ kg/m$^3$, $k_w = 1 \times 10^{-2}$ m/s. 

\begin{figure}[h!]
\centering
  \includegraphics[width=0.4\textwidth]{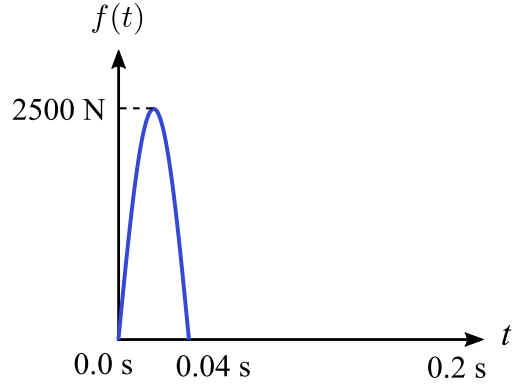}
\caption{Graphical depiction of the load profile used.}
\label{example2-load-profile}   
\end{figure}

As in \cite{pastor2000fractional, breuer1999quasi, markert2010comparison}, the problem domain is prescribed with null initial effective stress and water pressure. The dynamic load $f(t)$ is imposed on a strip footing on the top surface as shown in Figure \ref{example2-load-profile}. 
\begin{equation}
f(t) = 2500\sin(25\pi t) \mathbb{H},
\end{equation}
where $\mathbb{H}$ is equal to 1 if $t \leq 0.04$ s, and is zero if $t > 0.04$ s. The simulation time $t$ = 0.2 s and $\Delta t = 5\times 10^{-4}$ s. 

\begin{figure}[h!]
\centering
  \includegraphics[width=0.4\textwidth]{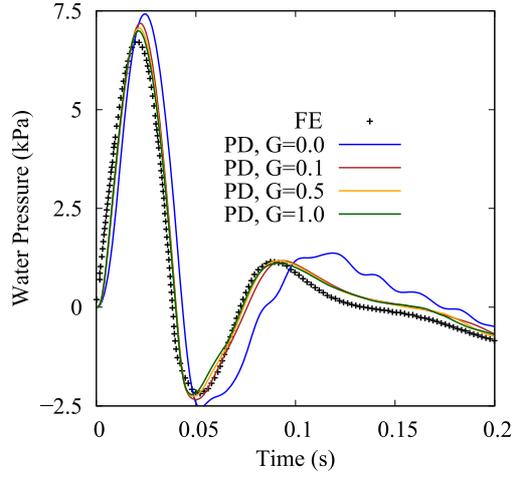}
\caption{Comparison of the PD solution and the FE solution of the water pressure at point A.}
\label{example2-results-pressure}   
\end{figure}

\begin{figure}[h!]
\centering
  \includegraphics[width=0.4\textwidth]{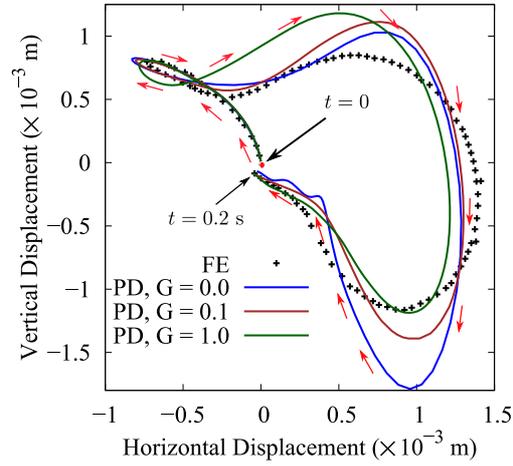}
\caption{Comparison of the PD solution the FE solution of the elliptical movement at point B.}
\label{example2-results-displacement}   
\end{figure}

Figure \ref{example2-results-pressure} compares the PD solutions with different values of $G$ and FE solutions \cite{markert2010comparison} of the water pressure at the point A. With $G$ = 0, the PD solution of the water pressure at point A slightly lags behind the FE solution and there are some oscillation at $t > 0.5$ s . The plot of water pressure with $G$ = 0.1 is smooth and closely matches the FE solution. The PD solutions with larger values of $G$ seem to generate mild changes in the plot of water pressure. Figure \ref{example2-results-displacement} compares the PD solutions with the FE solutions of the vertical displacement (heave) versus the horizontal displacement at point B. The results  show the expected elliptical motion at the point A under the dynamic load on the top boundary.



\begin{figure}[h!]
\centering
  \includegraphics[width=0.8\textwidth]{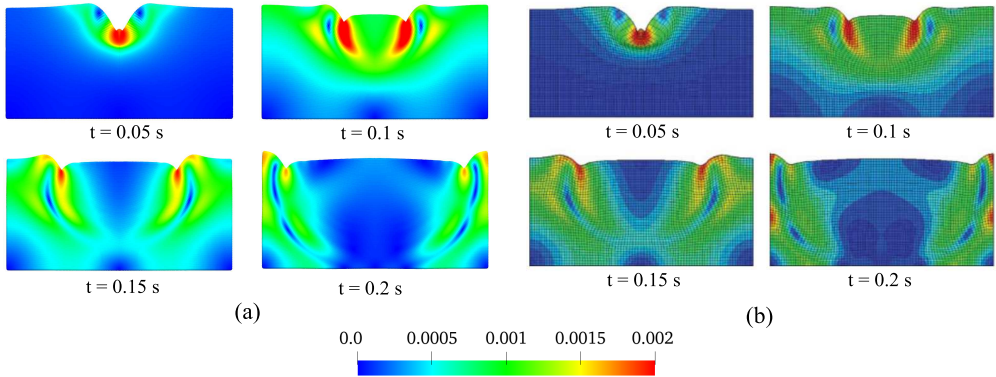}
\caption{Contours of the displacement magnitude (in m) and deformed configuration (magnification = $\times $500) predicted by PD (a) and FE (b) at $t = $ 0.05 s, 0.1 s, 0.15 s and 0.2 s.}
\label{example2-displacement-contours}   
\end{figure}

Figure \ref{example2-displacement-contours} plots the contour of the magnitude of displacement in the problem domain from the PD solution and the FEM solution at $t = $ 0.05 s, 0.1 s, 0.15 s and 0.2 s. In the PD solution, it is assumed $G$ = 0.05.

The results in Figure \ref{example2-displacement-contours} demonstrate that the PD solution is consistent with the FE solution in \cite{markert2010comparison}.  Both the PD and FE solutions show two-dimensional wave propagation through the bulk of the poroelastic medium.  with the surface showing elliptic motion by the dynamic load. The amplitude of this wave decreases as it travels farther away from the strip footing. Given the spike load profile, the deformation energy slowly dissipates as the wave moves through the problem domain that reduces the amplitude.

\subsection{Dynamic strain localization in a two-dimensional soil specimen}

\begin{figure}[h!]
\centering
  \includegraphics[width=0.25\textwidth]{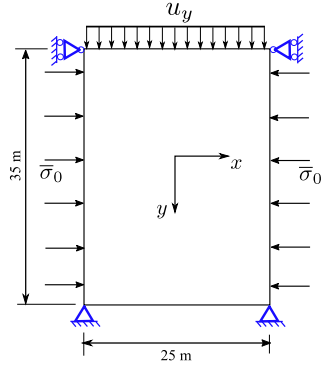}
\caption{Sketch of the geometry of the problem domain and boundary and loading conditions.}
\label{example3_problem_setup}
\end{figure}

In this example we simulate dynamic strain localization in a two-dimensional saturated soil sample subjected to a vertical compression. The problem geometry adopted and boundary conditions applied are shown in Figure \ref{example3_problem_setup}. 
The problem domain is discretized into 20,000 mixed material points with $\Delta x$ = 0.3 m and volume $0.027$ m$^3$. The material parameters used are : $\rho_s = 2000$ kg/m$^3$, $K = 2.5 \times 10^4$ kPa, $\mu_s = 1.154\times 10^4$ kPa, $\phi_0 = 0.3$,  initial pre-consolidation pressure $p_{c0} = -250.0$ kPa, $\kappa$ = 0.03, $\lambda$ = 0.10, $M$ = 1.0,  $\rho_w = 1000$ kg/m$^3$, $K_w = 2.0\times 10^5$ kPa, $k_w = 3\times 10^{-5}$ m/s. The horizon is set to 2.05$\Delta x$. The stabilization parameter $G$ = 0.025 is chosen to avoid excessive influence on the post-localization behavior \cite{silling2017stability}.

\begin{figure}[h!]
\centering
 \includegraphics[width=0.25\textwidth]{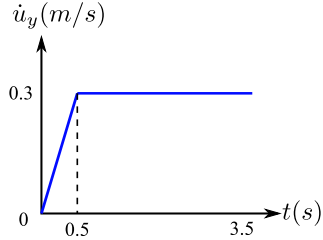}
\caption{Plot of the loading protocol for strain localization problem.}
\label{example3_loading_protocol}
\end{figure}

For the initial state the problem domain has an isotropic mean effective stress -100 kPa and zero water pressure. The skeleton boundary conditions are depicted in Figure \ref{example3_problem_setup}, where $\overline{\sigma}_0$ is a lateral confining pressure of 100 kPa. For the fluid phase, all boundaries are impermeable. The velocity load as shown in Figure \ref{example3_loading_protocol} is imposed on the top boundary. The simulation time at $t$ = 3.5 s and the time increment $\Delta t $ = 5$\times 10^{-3}$ s. 

\begin{figure}[h!]
\centering
 \includegraphics[width=0.5\textwidth]{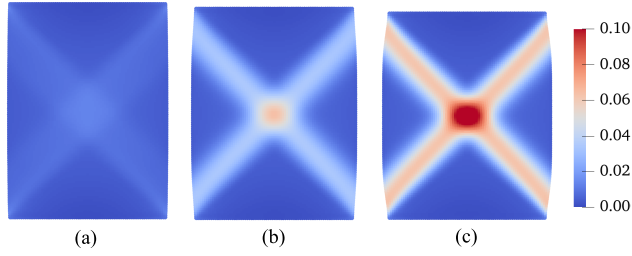}
\caption{Contours of equivalent shear strain ($\varepsilon_s$) at (a) $t$ = 1.0 s, (b) $t$ = 2.0 s, and (c) $t$ = 3.0 s (magnification = $\times$2). }
\label{example3-localization-contours-shear}
\end{figure}

\begin{figure}[h!]
\centering
 \includegraphics[width=0.5\textwidth]{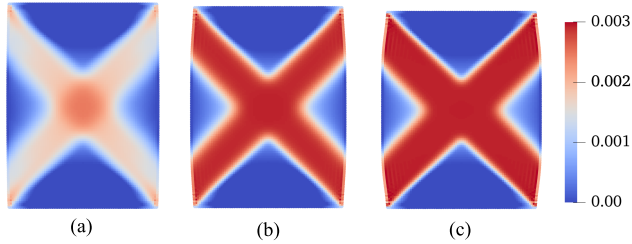}
\caption{Contours of plastic volume strain at (a) $t$ = 1.0 s, (b) $t$ = 2.0 s, and (c) $t$ = 3.0 s (magnification = $\times$2).}
\label{example3-localization-contours-volume}
\end{figure}

\begin{figure}[h!]
\centering
 \includegraphics[width=0.5\textwidth]{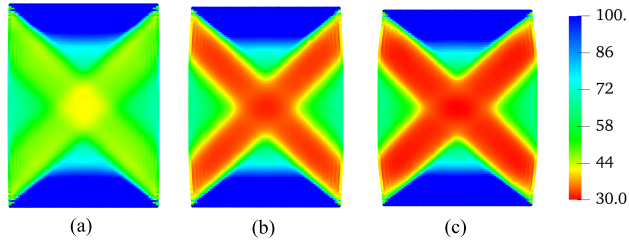}
\caption{Contours of water pressure (kPa) at (a) $t$ = 1.0 s, (b) $t$ = 2.0 s, and (c) $t$ = 3.0 s (magnification = $\times$2).}
\label{example3-localization-contours-pressure}
\end{figure}

Figures \ref{example3-localization-contours-shear}, \ref{example3-localization-contours-volume} and \ref{example3-localization-contours-pressure} draw the contours of the equivalent shear strain, plastic volume strain and water pressure at $t$ = 1.0 s,  2.0 s, and 3.0 s, respectively. The corresponding displacements on the top boundary are ${u}_y$ = 0.225 m, 0.52 m, 0.83 m, respectively. Here the equivalent shear strain is the second invariant of the strain tensor $\vec{\varepsilon}$, i.e., $\varepsilon_s $=$\sqrt{\frac{2}{3}}|\vec{\varepsilon} - \frac{1}{3}\mathrm{tr}(\vec{\varepsilon}) \vec{1}|$. 

The results show that the deformation and pressure have localized into symmetric banded zones. At $t$ = 1.0 s, plastic deformation has initiated and propagated diagonally through the specimen (Figures \ref{example3-localization-contours-shear} and \ref{example3-localization-contours-volume} (a)). We observe a similar behavior in the fluid pressure field (Figure \ref{example3-localization-contours-pressure}(a), where localization manifests as regions of decreased water pressure. In Figure \ref{example3-localization-contours-shear} (b) and (c) shows that plastic deformation progressively resolve into sharply defined zones of intense shear deformation. However, in the contours of plastic volume strain (see Figure \ref{example3-localization-contours-volume} (b) and (c)) and water pressure (see Figure \ref{example3-localization-contours-pressure}(b) and (c)) the banded zones appear to be more diffusive. The plastic volume change in the banded deformation is positive denoting dilatation. As such, the increase in skeleton volume leads to decrease in the water pressure inside the banded zone. The dilatation under dynamic loading can be expected for a moderately over-consolidated soil specimen. Due to the relatively large permeability chosen, the pore water can readily move into the zones of plastic dilatation, leading to more diffuse zones of water pressure. 

\subsubsection{Sensitivity to spatial discretization}

\begin{figure}[h!]
\centering
  \includegraphics[width=0.4\textwidth]{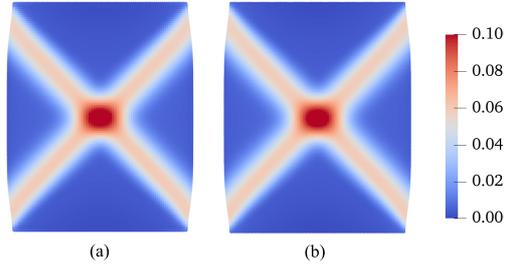}
\caption{Contours of equivalent shear strain ($\varepsilon_s$) at ${u}_y$ = 0.83 m for (a) coarse discretization and (b) fine discretization (magnification = $\times$2).}
\label{example3-mesh-contours-shear}
\end{figure}

\begin{figure}[h!]
\centering
  \includegraphics[width=0.4\textwidth]{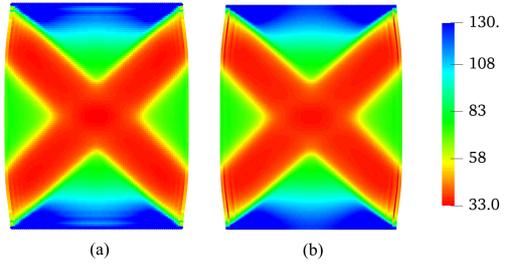}
\caption{Contours of water pressure (kPa) at ${u}_y$ = 0.83 m for (a) coarse discretization and (b) fine discretization (magnification = $\times$2).}
\label{example3-mesh-contours-pressure}
\end{figure}

In what follows, we present a discretization sensitivity analysis to demonstrate that the dynamic strain localization problem remains well-posed through the proposed nonlocal formulation. We rerun the numerical simulation with a fine spatial discretization. The fine spatial discretization consists of 43000 mixed material points with $\Delta x$ = 0.2 m.  For comparison, all material parameters and conditions remain the same. Figures \ref{example3-mesh-contours-shear} and \ref{example3-mesh-contours-pressure} compare the equivalent shear strain and water pressure respectively from the simulations with the coarse and fine discretizations at ${u}_y$ = 0.83 m. The results in Figures \ref{example3-mesh-contours-shear} show that the contours of the equivalent shear strain from both simulations are almost identical. The same conclusion can be drawn from Figure \ref{example3-mesh-contours-pressure} regarding the water pressure. It can be concluded that both he location and orientation of the shear band are independent of the spatial discretization.


\begin{figure}[h!]
\centering
 \includegraphics[width=0.42\textwidth]{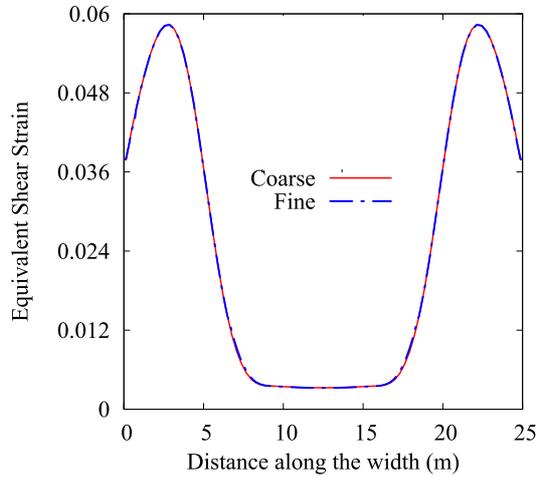}
\caption{Variations of equivalent plastic strain along a horizontal cross section at 10 m above the specimen center.}
\label{example3-mesh-plot-shear}
\end{figure}

\begin{figure}[h!]
\centering
 \includegraphics[width=0.4\textwidth]{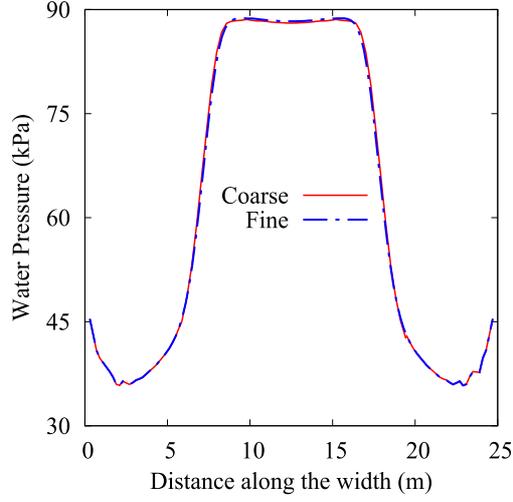}
\caption{Variations of water pressure along a horizontal cross section at 10 m above the specimen center.}
\label{example3-mesh-plot-pressure}
\end{figure}

In Figures \ref{example3-mesh-plot-shear} and \ref{example3-mesh-plot-pressure} we plot the variation of equivalent plastic shear strain and water pressure along a horizontal line at 10 m above the specimen center at $t$ = 3.0 s. The values of equivalent shear strain and water pressure are identical for both discretizations. 
It is known that in the dynamic strain localization analysis by FEM the mesh dependence of plastic strain and water pressure in the banded zone can be resolved by using a viscoplasticity model \cite{Loret1991Dynamic, shahbodagh2014dynamic}. However, the width of the banded zone still showed some sensitivity to the spatial discretization scheme (element size). 


\subsubsection{Influence of dynamic loading}

\begin{figure}[h!]
\centering
  \includegraphics[width=0.5\textwidth]{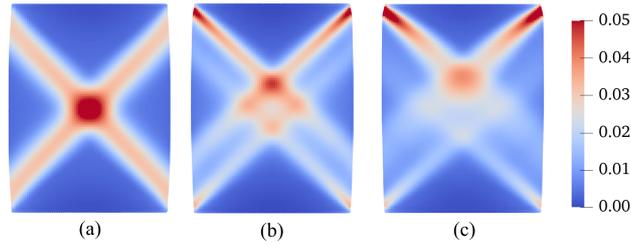}
\caption{Contours of equivalent shear strain ($\varepsilon_s$) at ${u}_y$ = 0.5 m from simulations with three loading rates: (a) $\dot{u}_y$ = 0.3 m/s, (b) $\dot{u}_y$ = 0.9 m/s, and (c) $\dot{u}_y$ = 1.5 m/s (magnification = $\times$2). }
\label{example3-loadrates-contours-shear}
\end{figure}

\begin{figure}[h!]
\centering
  \includegraphics[width=0.5\textwidth]{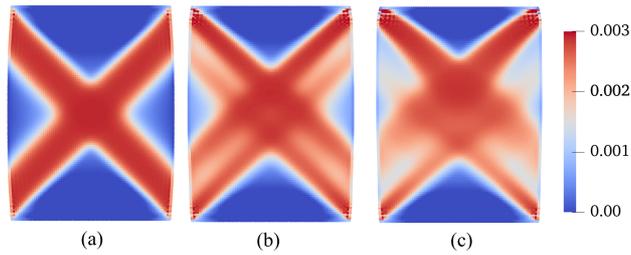}
\caption{Contours of plastic volume strain at ${u}_y$ = 0.5 m from simulations with three loading rates: (a) $\dot{u}_y$ = 0.3 m/s, (b) $\dot{u}_y$ = 0.9 m/s, and (c) $\dot{u}_y$ = 1.5 m/s (magnification = $\times$2).}
\label{example3-loadrates-contours-volume}
\end{figure}

\begin{figure}[h!]
\centering
  \includegraphics[width=0.5\textwidth]{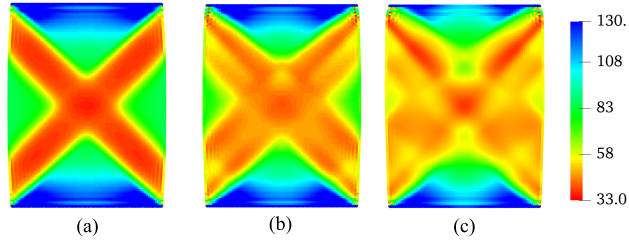}
\caption{Contours of water pressure (kPa) at ${u}_y$ = 0.5 m from simulations with three loading rates: (a) $\dot{u}_y$ = 0.3 m/s, (b) $\dot{u}_y$ = 0.9 m/s, and (c) $\dot{u}_y$ = 1.5 m/s (magnification = $\times$2).}
\label{example3-loadrates-contours-pressure}
\end{figure}

We investigate the influence of dynamic loading rates on the coupled response during dynamic strain localization. The base simulation with $\dot{u}_y = 0.3$ m/s is repeated with $\dot{u}_y = 0.9$ m/s and $\dot{u}_y$ = 1.5 m/s respectively. All other material parameters and conditions remain the same. The results are compared at an identical displacement of the top boundary. Figures \ref{example3-loadrates-contours-shear}, \ref{example3-loadrates-contours-volume} and \ref{example3-loadrates-contours-pressure} show the contours of the equivalent shear strain, plastic volume strain and water pressure respectively at ${u}_y$ = 0.5 m on the top boundary for three different dynamic loading rates.  

The results these Figures show the loading rate impact the formation of banded deformation and pressure zones. For $\dot{u}_y$ = 0.3 m/s, it is apparent that the contour of plastic deformation as shown in Figures \ref{example3-loadrates-contours-shear} and \ref{example3-loadrates-contours-volume} (a) and the contour of water pressure  as shown in Figure \ref{example3-loadrates-contours-pressure}(a) have localized into a single pair of two conjugate banded zones. For the simulations with two larger loading rates, $\dot{u}_y$ = 0.9 m/s and 1.5 m/s, two distinct pairs of localized deformation bands are formed in the solid skeleton as shown in Figures \ref{example3-loadrates-contours-shear} and \ref{example3-loadrates-contours-volume} (b). Consistent with the skeleton deformation, the contour of water pressure also shows two pairs of banded zones as shown in Figure \ref{example3-loadrates-contours-pressure}(b). It can be concluded from the results in Figures \ref{example3-loadrates-contours-shear}, \ref{example3-loadrates-contours-volume} and \ref{example3-loadrates-contours-pressure} that the dynamic loading rate impacts the number, location and orientation of shear bands in unsaturated porous media.

\begin{figure}[h!]
\centering
  \includegraphics[width=0.4\textwidth]{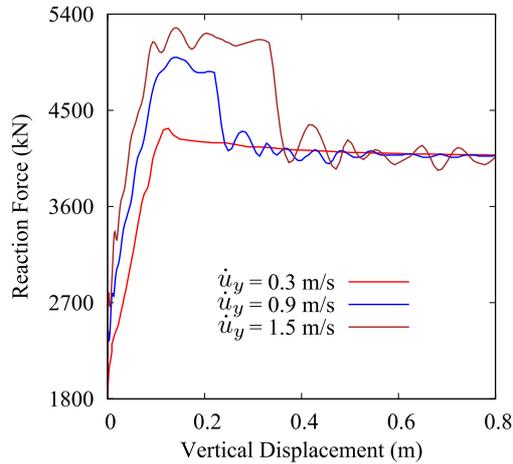}
\caption{Plot of reaction force over the vertical displacement at the top boundary for simulations with three loading rates.}
\label{example3-loadrates-reaction}
\end{figure}

Figure \ref{example3-loadrates-reaction} plots the reaction force over the applied vertical displacement for the three different loading rates. The results show that the specimen under a higher dynamic loading rate shows a higher peak value in the loading capacity. The oscillations in the reaction force curves with larger loading rates may be correlated to locally undrained condition under these loading condition and thus less dilatation in the specimen. 

\section{Closure}

In this article we propose a stabilized computational nonlocal poromechanics model for dynamic strain localization in saturated porous media. The stabilized coupled nonlocal model is solved using a Lagrangian-Eulerian meshless method with an implicit time integration scheme. Parallel computing is adopted for computing efficiency. As a new contribution, we present a theoretical proof of the zero-energy modes associated with the multiphase correspondence principle. We propose a remedy based on the energy method to circumvent zero-energy modes in the solid deformation and fluid flow. We present a method to determine the stabilization parameter $G$ for both the solid deformation and fluid flow process. The stabilized coupled nonlocal model for saturated porous media under dynamic loading can be readily extended to model dynamic problems in unsaturated porous media. We have validated the coupled stabilized nonlocal formulation by comparing numerical results with analytical and finite element solutions for dynamic problems in saturated porous media. Numerical examples are conducted to demonstrate the robustness of the coupled nonlocal dynamic model for dynamic strain localization analysis of saturated porous media. 

\section*{Acknowledgments}

The work presented in this article has been supported by the US National Science Foundation under contract numbers 1659932 and 1944009. 

\bibliography{pd-stabilization}

\end{document}